\documentclass{article}

\usepackage{arxiv}

\usepackage[utf8]{inputenc} 
\usepackage[T1]{fontenc}    
\usepackage{hyperref}       
\usepackage{url}            
\usepackage{booktabs}       
\usepackage{amsfonts}       
\usepackage{nicefrac}       
\usepackage{microtype}      
\usepackage{lipsum}
\usepackage{natbib}
\usepackage{amssymb}
\usepackage{amsmath}
\usepackage{amscd}
\usepackage{dsfont}
\usepackage{graphicx}
\usepackage{lipsum}
\usepackage{enumitem}
\graphicspath{{./Figures_final/}} 
\usepackage[noadjust]{cite}
\usepackage{amsthm}
\usepackage{caption}
\usepackage{multirow}
\usepackage[toc,page,header]{appendix}
\usepackage{array}
\usepackage{pifont}
\usepackage[linesnumbered,ruled,vlined]{algorithm2e}
\usepackage{hyperref}
\usepackage{color}
\usepackage{indentfirst}
\usepackage{optidef}

\hypersetup{colorlinks,linkcolor={blue},citecolor={blue},urlcolor={red}}
\newcommand{\bs}{\boldsymbol}

\DeclareMathOperator*{\argmin}{arg\,min}

\SetCommentSty{mycommfont}

\SetKwInput{KwInput}{Input}                
\SetKwInput{KwOutput}{Output}              

\renewcommand{\thealgocf}{}

\renewcommand{\u}[0]{\bold{u}}                        
\renewcommand{\v}[0]{\bold{v}}                        
\renewcommand{\b}[0]{\bold{b}}                        
\renewcommand{\d}[0]{\bold{d}}                        
\renewcommand{\P}[0]{\bold{P}}                       
\renewcommand{\next}[1]{#1_{k}}                         

\newcommand{\m}[0]{\bold{m}}                        
\newcommand{\A}[0]{\bold{A}}                        
\newcommand{\Amu}{\A(\m)\u} 
\newcommand{\Pu}{\P\u} 
\newcommand{\diag}[0]{\text{diag}}                  
\renewcommand{\c}[0]{\bold{c}}                        
\DeclareMathOperator*{\minimize}{minimize}

\begin{document}
\title{Robust elastic full waveform inversion using alternating direction method of multipliers with reconstructed wavefields}

\author{
 Kamal Aghazade \\
 Formerly at the Institute of Geophysics, University of Tehran, Tehran, Iran \\
 Currently at the Institute of Geophysics, Polish Academy of Sciences, Warsaw, Poland. \\
 \texttt{aghazade.kamal@igf.edu.pl} \\
   \And
 Ali Gholami \\
  Institute of Geophysics, Polish Academy of Sciences, Warsaw, Poland, \\
  \texttt{agholami@igf.edu.pl} (corresponding author)\\
  \And
 Hossein S. Aghamiry \\
Charité Universitätsmedizin Berlin \\
Center for Biomedicine, Berlin, Germany.
  \texttt{hossein.aghamiry@charite.de} \\
  \AND
  Hamid Reza Siahkoohi \\
    Institute of Geophysics, University of Tehran, Tehran, Iran. \\
   \texttt{hamid@ut.ac.ir}
}

\maketitle
\begin{abstract}
Elastic full-waveform inversion (EFWI) is a process used
to estimate subsurface properties by fitting seismic data while
satisfying wave propagation physics. The problem is formulated as a least-squares data fitting minimization problem with
two sets of constraints: Partial-differential equation (PDE)
constraints governing elastic wave propagation and physical
model constraints implementing prior information. The alternating direction method of multipliers is used to solve the
problem, resulting in an iterative algorithm with well-conditioned subproblems. Although wavefield reconstruction is the
most challenging part of the iteration, sparse linear algebra
techniques can be used for moderate-sized problems and frequency domain formulations. The Hessian matrix is blocky
with diagonal blocks, making model updates fast. Gradient
ascent is used to update Lagrange multipliers by summing
PDE violations. Various numerical examples are used to investigate algorithmic components, including model parameterizations, physical model constraints, the role of the Hessian
matrix in suppressing interparameter cross-talk, computational efficiency with the source sketching method, and the
effect of noise and near-surface effects.  
\end{abstract}
%
%
\section{Introduction}
Full waveform inversion (FWI) is a widely used seismic imaging technique for estimating the elastic properties of the subsurface by inverting seismic data recorded at or near the surface. 
The modeling process is based on solving the wave equation \citep{Tarantola_1984_ISR}. 
The inverse process involves a non-linear optimization problem, and due to the computational cost, it typically uses a local optimization strategy to iteratively modify the properties of the medium \citep [for a comprehensive review, see ][] {Virieux_2009_OFW}.
The fidelity of the resulting property estimates hinges upon the interplay between accurate forward modeling, which meticulously incorporates the physics of wave propagation, and the efficacy of the chosen inverse optimization method.

Initially, FWI applications considered either the acoustic approximation \citep{Tarantola_1984_ISR,Gauthier_1986_TNI} or the elastic isotropic approximation \citep{Mora_1988_EWI} of the Earth's interior. Acoustic FWI assuming a constant density was preferred for field records since the unconverted P-wave velocity ($\bold{V}_\text{P}$) could describe most features in the data \citep[e.g., ][]{Ravaut_2004_MSI,Operto_2006_CIM}.
However, traditional acoustic FWI methods overlook the presence of S-waves in the data \citep{Prieux_2011_FAI}. This limitation leads to inaccurate modeling of physics and introduces residual data artifacts \citep{Mulder_2008_ESI}. Previous studies have attempted to address these issues by considering density as a proxy variable in acoustic FWI, which partially accounts for dynamic elastic effects and helps mitigate overfitting of elastic data caused by incorrect $\bold{V}_\text{P}$ variations \citep{borisov2014acoustic}.
Despite these efforts, acoustic inversion still faces challenges in accurately modeling amplitude variation with offset effects, particularly for wide aperture data  \citep{Barnes_2008_FWI,Barnes_2009_DAA}. This is due to amplitude errors arising from the differences between the recorded particle velocities and pressure fields, as well as the directivity of the sources and receivers. 
Addressing these limitations and achieving a more accurate subsurface characterization require the development of robust and efficient inversion techniques that account for the complete elastic forward modeling.
%
Recent advances in computational power and the acquisition of multicomponent data have enabled multiparameter inversion, where different classes of parameters can be incorporated into the FWI procedure \citep{Sears_2008_EFW, Prieux_2013_MFWa,Prieux_2013_MFWb,Operto_2013_TLE,vigh2014elastic}.

Despite increasing the accuracy of physics of wave propagation, transitioning from acoustic to elastic inversion poses challenges that must be addressed appropriately. The sensitivity of the inversion for one parameter class may differ significantly from another, and adding more parameters introduces more degrees of freedom, potentially increasing the degree of ill-posedness. Cross-talk, which refers to the coupling between different parameter classes as a function of the scattering angle, is another issue in multiparameter inversion. Parameters with different natures may have similar signatures on the data, and this can be quantified using tools such as radiation pattern analysis \citep{Forgues_1997_PSA,kazei2019scattering} or sensitivity kernel analysis \citep{pan2019interparameter}.
Several strategies have been proposed to address parameter cross-talk, including model/data-driven approaches, model reparameterization, the use of the inverse of the Hessian matrix, {and decomposition of the P-and S-wave modes \citep{Gholami_2013_WPA1,Metivier_2014_MFW,wang2017elastic,wang2018elastic,wang2019elastic}}.

Once the forward modeling operator is established, the parameter estimation task transforms into a nonlinearly constrained optimization problem \citep{haber2000optimization}. This optimization objective aims to minimize the data misfit while adhering to nonlinear constraints dictated by the wave equation. Various techniques, including Lagrangian, penalty, and Augmented Lagrangian (AL) methods, can be employed to address such constrained problems \citep{nocedal2006numerical}. Additionally, optimization can be approached through either reduced space or full space methods. In the former, optimization exclusively involves model parameters, while the latter includes wavefields and Lagrange multipliers. Owing to superior memory efficiency, the majority of FWI algorithms opt for the reduced space approach \citep{bunks1995multiscale,Operto_2006_CIM,brossier2009seismic,kohn2012influence,duan2016elastic}.

However, reduced space FWI encounters cycle skipping (local minima), a challenge exacerbated in Elastic FWI (EFWI) by the propagation of short-wavelength S-waves. To address this issue, early studies primarily relied on low-frequency data (below 3 Hz) \citep[e.g.,][]{choi2008frequency, brossier2009seismic, kohn2012influence} or implemented model-driven workflows. These workflows, such as using hydrophone data to recover long-to-intermediate wavelengths of S-wave velocity ($\bold{V}_\text{S}$), were instrumental in establishing optimal initial models for subsequent inversion phases \citep[see e.g.,][]{Prieux_2013_MFWb}. Inversion processes often adopt a multiscale approach, commencing with low frequencies and gradually progressing to higher frequencies, to mitigate the impact of cycle skipping \citep{bunks1995multiscale}.

When confronting seismic data lacking low-frequency components, EFWI may utilize the low-frequency component of the damped wavefield in the Laplace-Fourier domain \citep{jun2014laplace}. Nevertheless, the presence of undesirable cross-correlation components can distort the gradient, particularly in scenarios like low-velocity zones beneath salt, potentially yielding low-resolution and inaccurate models \citep{kwon2017waveform}. Machine- learning (ML) tools have been explored to extrapolate low-frequency data \citep{sun2021deep}; however, their utilization often involves significant computational costs and remains an active area of research. Phase tracking, as used  by \citet{li2015phase}, enhances inversion robustness in the absence of low frequencies, particularly in dispersion-free media. To address cycle skipping, researchers, such as \citet{chi2014full} and \citet{zhang2022elastic}, have turned to envelope-based inversion through wave mode decomposition. However, envelope inversion grapples with instabilities in both the gradient term and the adjoint source \citep{xiong2023improved}. Combining envelope inversion with deconvolution, as proposed by \citet{chen2022elastic}, considers phase shifts and mitigates artificial side lobes in data due to the loss of low-frequency components. Phase correction methods, employed by researchers like \citet{chi2014full} and \citet{hu2022frequency}, address phase differences between modeled and synthetic data. Another strategy involves modifying the misfit function; for instance, \citet{zhang2019local} introduced a local cross-correlation objective function to overcome the susceptibility to cross-talk from neighboring arrivals. The application of optimal transport misfit functions, proposed in some studies \citep{marty2022elastic,borisov2022graph}, faces challenges, especially when comparing signed and oscillatory signals in elastic media.

The challenge of local minima in reduced space FWI can also be effectively addressed using penalty or AL formulations in the full space. For acoustic FWI, \citet{VanLeeuwen_2013_MLM} applied the penalty formulation, leading to the Wavefield Reconstruction Inversion (WRI) algorithm, whereas  \citet{aghamiry2019improving, aghamiry2020multiparameter} utilized the AL formulation, resulting in the iteratively refined WRI (IR-WRI) algorithm. Extensive numerical results have demonstrated a substantial enhancement in stability and convergence with these algorithms compared to traditional reduced-space FWI \citep[see][ and references therein]{Operto_2023_FWI}. Moreover, algorithms based on AL methods outperform penalty methods in terms of stability and convergence \citep[for a detailed comparison between them, see][ ]{Gholami_2023_FWI}.

The AL (AL) formulation provides a comprehensive framework for addressing challenges in EFWI, encompassing issues of ill-posedness, cross-talk, and convergence. Leveraging the Alternating Direction Method of Multipliers (ADMM) \citep{Gabay_1976_ADA, boyd2011distributed}, the resulting algorithm facilitates the integration of the Hessian matrix to mitigate cross-talk, incorporates regularization techniques to alleviate ill-posedness \citep{lin2014acoustic}, and demonstrates robust convergence even from inaccurate initial models, all while maintaining computational efficiency. To adapt the ADMM method to elastic media, we explore various parameterizations, including Lamé parameters and velocities, with density set as a constant. As demonstrated later, in both cases, the associated Gauss-Newton Hessian matrix exhibits a block-diagonal sparse structure that allows explicit inversion, enabling an exploration of its role in reducing interparameter cross-talk. Additionally, we introduce a source sketching method, previously developed for acoustic FWI, to enhance the computational efficiency of the algorithm.

Furthermore, we incorporate physical constraints into the EFWI inversion process to promote physically plausible models \citep{duan2016elastic}. One such methodology involves the use of seismic facies information built by clustering seismic characteristics and spatial coherence within the framework of FWI \citep{zhang2018multiparameter}. Another approach involves using the relationships between the parameters in the elastic medium based on borehole data or empirical equations derived from laboratory measurements \citep{brocher2005empirical}. We optimize the model parameters so that they are at the intersection of these constraints, and this problem is efficiently addressed using the ADMM method. We present numerical examples using benchmark models to demonstrate the performance of the proposed EFWI algorithm.

%
%
\section{Preliminaries}
The basic definitions and symbols used in this paper are as follows. The field of real and complex numbers are denoted by $\mathbb{R}$ and $\mathbb{C}$, respectively. Vectors and matrices are represented by bold lowercase and uppercase letters, respectively.
The number of discrete model parameters (for each class), the number of receivers, and the number of sources are denoted by $n$, $n_r$, and $n_s$, respectively. The angular frequency is denoted by $\omega = 2\pi f$, where $f$ is the frequency.
The identity and diagonal matrices are denoted by $\bold{I}$ and $\text{diag}(\bullet)$, respectively. The second-order partial derivative with respect to variable $i$ is represented by $\partial_{ii}$, and the mixed derivative with respect to variables $i$ and $j$ is represented by $\partial_{ij}$. The conjugate transpose of a matrix/vector is denoted by the superscript $T$. The symbol $\circ$ denotes the Hadamard (element-wise) multiplication. In this paper, $(\bullet)^k$ represents the value of $(\bullet)$ at iteration $k$. 
%
%
%
%
\section{Theory}
We consider the following frequency-domain isotropic elastic wave equation:
\begin{subequations}\label{El_forward}
\begin{align}
\bs{\rho} \omega^{2}\bold{u}_x+(\bs{\lambda}+2\bs{\mu})\partial_{xx}\bold{u}_x+\bs{\mu} \partial_{zz}\bold{u}_x+(\bs{\lambda}+\bs{\mu})\partial_{xz}\bold{u}_z = \bold{b}_x, \label{EL_bx} \\
\bs{\rho} \omega^{2}\bold{u}_z+(\bs{\lambda}+2\boldsymbol{\mu})\partial_{zz}\bold{u}_z+\bs{\mu} \partial_{xx}\bold{u}_z+(\bs{\lambda+\mu})\partial_{xz}\bold{u}_x = \bold{b}_z, \label{EL_bz}
\end{align}
\end{subequations}
where $\bs{\rho} \in \mathbb{R}^{n \times 1}$ is mass density, $\bs{\lambda} \in \mathbb{R}^{n \times 1}$ and $\bs{\mu} \in \mathbb{R}^{n \times 1}$ denote Lamé parameters, $\bold{u}_x \in \mathbb{C}^{n \times 1}$ and $\bold{u}_z \in \mathbb{C}^{n \times 1}$ are horizontal and vertical particle displacements, and $\bold{b}_x \in \mathbb{C}^{n \times 1}$, $\bold{b}_z \in \mathbb{C}^{n \times 1}$ are the source terms. Equation \ref{El_forward} can be written in compact algebraic form as:
\begin{equation}\label{Aub}
\small
\Amu = \b, 
\end{equation} 
where 
\begin{equation}
 \m = 
\begin{pmatrix} 
 \bs{\lambda}\\
 \bs{\mu}\\
  \bs{\rho}
\end{pmatrix},
\quad
\bold{u} = 
\begin{pmatrix} 
\u_x\\
\u_z
\end{pmatrix},
\quad
\b = 
\begin{pmatrix} 
\b_x\\
\b_z
\end{pmatrix},
\end{equation}
and
\begin{equation}
\noindent 
\bold{A}( \m) = \begin{pmatrix}
\omega^2\diag(\bs{\rho})+\diag(\bs{\varpi})\partial_{xx}+\diag(\bs{\mu})\partial_{zz} ~~~~ \diag(\bs{\lambda}+\bs{\mu})\partial_{xz} \\
\diag(\bs{\lambda}+\bs{\mu})\partial_{xz} ~~~~ \omega^2\diag(\bs{\rho})+\diag(\bs{\varpi})\partial_{zz}+\diag(\bs{\mu})\partial_{xx} 
\end{pmatrix},
\end{equation}
denotes the complex-valued impedance matrix, where $\bs{\varpi} = \bs{\lambda}+2\bs{\mu}$.

\begin{figure}[!h] 
\centering 
\includegraphics[scale=1]{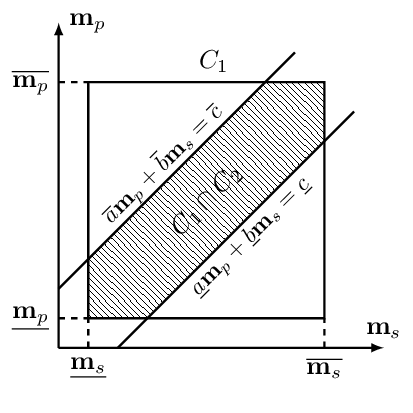}
\caption{The sets $C_1$ and $C_2$ and the related parameters. The filled area shows the desired set.}
\label{Fig:set}
\end{figure}

In order to estimate the subsurface model parameters $\mathbf{m}$, we formulate the elastic FWI as the following constrained optimization problem:
\begin{equation}\label{const_fwi_cnvx}
\begin{split}
&\minimize_{\m,\u}~\frac{1}{2}\|\P \u-\d\|_{2}^2 \\
& \text{subject to}\quad \Amu = \b, \\
& \hspace{1.7cm} ~\m \in C_1 \cap C_2
\end{split}
\end{equation}
where
\begin{equation}
\P = 
\begin{pmatrix} 
\tilde{\P}& \bold{0}\\
\bold{0}& \tilde{\P}
\end{pmatrix},
\quad
\d = 
\begin{pmatrix} 
\d_x\\
\d_z
\end{pmatrix},
\end{equation}
and the observation operator $\tilde{\mathbf{P}}$ samples the wavefield at receiver locations.
$C_1$ and $C_2$ are two convex sets built from prior physical information about the model parameters. Specifically, $C_1$ defines a box in the plane formed by ($\mathbf{m}_p$,$\mathbf{m}_s$) with lower bound $(\underline{\mathbf{m}_p},\underline{\mathbf{m}_s})$ and upper bound $(\overline{\mathbf{m}_p},\overline{\mathbf{m}_s})$. Also, $C_2$ is defined as the set of points between two lines $\underline{a} \mathbf{m}_p+\underline{b}\mathbf{m}_s=\underline{c}$ and $\overline{a} \mathbf{m}_p+\overline{b}\mathbf{m}_s=\overline{c}$, where the parameters $(\underline{a},\underline{b},\underline{c})$ and $(\overline{a},\overline{b},\overline{c})$ can be found according to the borehole information or empirical relations \citep{brocher2005empirical,duan2016elastic}. The intersection of $C_1$ and $C_2$ forms a convex set as shown in  Figure~\ref{Fig:set}.

Satisfying such constraints can be achieved by solving the corresponding convex feasibility problem by projecting the updated models onto the intersection of nonempty, closed, and convex sets. One simple approach to solve the feasibility problem is von Neumann's alternating projection method (also known as the projection onto convex sets, POCS), which was used by \citet{baumstein2013pocs} for the case of anisotropic acoustic inversion. Another variant of alternating projection with slight differences by introducing auxiliary variables is Dykstra's algorithm \citep{dykstra1983algorithm}.

From equation \ref{const_fwi_cnvx}  one can derive the AL \citep{Powell_1969_NLC}
\begin{equation}\label{ALF}
\small
\mathcal{L}_{\beta}(\u,\m,\bs{\eta}) = \frac{1}{2}\|\Pu-\bold{d}\|_{2}^2-\boldsymbol{\eta}^{T}\left[\Amu - \b\right]+\frac{\beta}{2} \|\Amu- \b\|_{2}^{2}.
\end{equation} 
In this equation, $\mathcal{L}_{\beta}(\u,\m,\bs{\eta})$ consists of three terms. The first term is the data misfit term, quantifying the discrepancy between the modeled data $\Pu$ and the observed data $\bold{d}$. The second term, a Lagrangian term, introduces the Lagrange multiplier vector $\bs{\eta}$ to enforce satisfaction of the wave equation constraint. To enhance optimization stability, the third term, a penalty term, is added, controlled by the penalty parameter $\beta>0$. The adjustment of $\beta$ penalizes the Lagrange multipliers, contributing to optimization stability \citep{Gholami_2023_FWI}.\\
Notably, the choice of the negative sign in the Lagrangian term aligns with the convention set by \citet{Powell_1969_NLC}. Although the optimization of the Lagrangian function alone might be prone to instability, the introduction of the penalty term through the AL technique significantly improves convergence. Importantly, this technique surpasses the penalty method, as the estimates of Lagrange multipliers in the AL approach tend to be more accurate. It is crucial to highlight that setting the Lagrange multiplier term to zero reverts to the penalty formulation used in the WRI approach \citep{VanLeeuwen_2013_MLM}. The optimization of the AL function is effectively addressed through dual ascent methods within the framework of the wavefield-oriented ADMM \citep{boyd2011distributed,Gholami_2023_FWI}:
\begin{subequations}\label{ADMM}
\begin{align}
\small
\u^{k+1} &= \argmin_{\u} ~\mathcal{L}_{\beta}(\u,\m^{k},\bold{s}^{k}), \label{ADMM_u} \\
\m^{k+1} &= \argmin_{\m \in C_1 \cap C_2}~ \mathcal{L}_{\beta}(\u^{k+1},\m,\bold{s}^{k}), \label{ADMM_m} \\
\bold{s}^{k+1} &= \bold{s}^{k}+\b-\A(\m^{k+1})\u^{k+1}, \label{ADMM_eta}
\end{align}
\end{subequations}  
where $\bold{s}^{k}=\frac{1}{\beta}\bs{\eta}^{k}$ is the scaled dual variable.
In the following subsections, we provide the detailed analysis for solving these subproblems.
%


\subsection{Wavefield reconstruction step}
The optimization problem in \ref{ADMM_u} is quadratic in $\bold{u}$, and its minimization leads to the following linear system of equations:
\begin{equation}\label{u_sub}
\small
\left[ \beta \A^{T}(\m^{k})\A(\m^{k}) + \P^{T}\P \right] \u^{k+1} = \beta \A^{T}(\m^{k})\left(\b+\bold{s}^{k}\right)+\P^{T}\d.  
\end{equation} 
where the size of this system is twice the size of the system that appears in acoustic inversion \citep{VanLeeuwen_2013_MLM, aghamiry2019improving}. If the problem size is small to moderate, this system can be solved explicitly using LU factorization. However, for larger problems, iterative methods with appropriate preconditioners should be used to solve this data-assimilated wave equation. 
For 3D applications, the time-domain implementation is recommended. However, tackling equation \ref{u_sub} in the time-domain poses a substantial challenge due to the absence of explicit time-stepping and the escalating computational and memory demands associated with larger models. This challenge arises from the inherently high-dimensional nature of the wavefield. A more efficient alternative approach involves adopting an equivalent multiplier-oriented formulation. In this strategy, we first solve a normal system of data size to obtain least-square multipliers. Subsequently, the wavefields are computed straightforwardly using time-stepping \citep[see][]{gholami2022extended, Gholami_2023_FWI}.

%
%

%
\subsection{Model estimation step}
\begin{figure}[!h]
	\centering 
	\includegraphics[scale=.6]{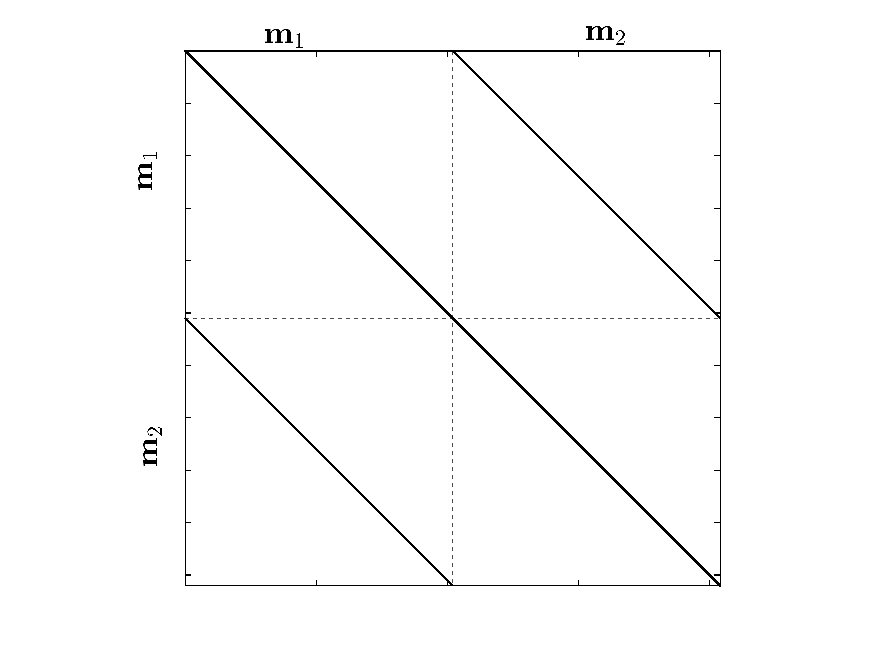}
	\caption{Hessian structure for two-parameters problem. The main-diagonal blocks and off-diagonal blocks of  
		$\beta \bold{L}^{T}\bold{L}$ respectively include 
		$\frac{\partial^2 \mathcal{L}}{\partial \bold{m}_{i} \partial \bold{m}_{i}}$
	and
		$\frac{\partial^2 \mathcal{L}}{\partial \bold{m}_{i} \partial \bold{m}_{j}}$, in which 
		{ $i,j$} denote the index of each parameter class.
	}
	\label{Fig:Hessian_stuct}
\end{figure}

Different parameterizations can be used to update model parameters in elastic waveform inversion \citep{Tarantola_1986_SNL}. The Lamé parameters ($\bs{\lambda},\bs{\mu}$) and density ($\bs{\rho}$) are commonly used for defining an elastic isotropic medium. However, density reconstruction poses challenges due to parameter trade-off, poor sensitivity of waveform mismatch to the density perturbation, and difficulties in retrieving the low-wavenumber component of density \citep{kohn2012influence,sun2017density,pan2018elastic}. 
In this study, we consider a constant density medium and explore parameterizations based on ($\bs{\lambda},\bs{\mu}$) and ($\bold{V}_{P}^2,\bold{V}_{S}^2$) for model update, where for the $\m = (\bold{V}_{P}^2,\bold{V}_{S}^2)$ parameterization, the relations $\bs{\lambda} = \bs{\rho}\bold{V}_{P}^2-2\bold{V}_{S}^2$ and $\bs{\mu} = \bs{\rho}\bold{V}_{S}^2$
are used. 
Regarding the model estimation problem, we can write:
  \begin{equation}\label{L_lame}
  	\bold{A}(\bold{m})\bold{u}^{k+1}-\bold{b}-{\bold{s}}^{k} = \underbrace{\begin{bmatrix}
  			\bold{L}_{11} &  \bold{L}_{12}\\
  			\bold{L}_{21} & \bold{L}_{22}
  	\end{bmatrix}}_{\bold{L(\bold{u}^{k+1})}}\underbrace{\begin{bmatrix}
  		\bold{m}_1 \\ \bold{\bold{m}}_2
  	\end{bmatrix}}_{\bold{m}}-\underbrace{\begin{bmatrix}
  		\bold{y}_x\\
  			\bold{y}_z\\
  	\end{bmatrix}}_{\bold{y}(\bold{u}^{k+1})},
  \end{equation}
  Then, for each mentioned parameterization, the updated parameters are obtained by solving the following least squares problem:
\begin{equation}\label{m_update_linearOp}
\small
\bold{m}^{k+1} = \argmin_{\bold{m}}~ \|\bold{L}\bold{m}-\bold{y} \|_{2}^{2},
\end{equation}
where $\bold{L} \equiv \bold{L}(\bold{u}^{k+1})$ and $\bold{y} \equiv \bold{y}(\u^{k+1})$. This least-squares problem admits the closed-form solution:
\begin{equation}\label{m_update_closed}
\small
\m^{k+1} = \left[ \bold{L}^{T}\bold{L} \right] ^{-1}\bold{L}^{T} \bold{y}.
\end{equation}
The structure of the operator $\bold{L}$ and vector $\bold{y}$ are presented in Table~\ref{Table:L_lame_rho} for different parameterizations. The operator $\bold{L}$ is a 2 by 2 block matrix with diagonal blocks. Consequently, the sparsity pattern and straightforward calculation process facilitate the explicit computation of the inverse of the Gauss-Newton Hessian matrix $\bold{L}^{T}\bold{L}$ (Figure~\ref{Fig:Hessian_stuct}).


We also consider the parameterization $(\bold{V}_{P},\bold{V}_{S})$. Mathematical formulas for the model update in this case are presented in Appendix A. 
In Appendix B, we provide the detailed formulation for solving the model subproblem by including physical constraints.
\begin{table*}[!h]
\centering
\caption {The structure of matrix $\bold{L}$ and vector $\bold{y}$ for ($\bs{\lambda}$, $\bs{\mu}$) and ($\bold{V}_{P}^{2},\bold{V}_{S}^{2}$) parameterization. }
\footnotesize
\begin{tabular}{|c|c|c|}
\hline
 &   $\bold{L}(\u)\in \mathbb{C}^{2n\times 2n}$ &  $\bold{y}(\u)\in \mathbb{C}^{2n\times 1}$ 
\\ \hline
$(\bs{\lambda},\bs{\mu})$ 
& 
$\begin{bmatrix}
\diag(\partial_{xx}\bold{u}_{x})+\diag(\partial_{xz}\bold{u}_z) &  \diag(2\partial_{xx}\bold{u}_{x})+\diag(\partial_{zz}\bold{u}_{x})+\diag(\partial_{xz}\bold{u}_z) \\
 \diag(\partial_{zz}\bold{u}_z)+\diag(\partial_{xz}\bold{u}_{x}) & \diag(2\partial_{zz}\bold{u}_{z})+\diag(\partial_{xx}\bold{u}_z)+\diag(\partial_{xz}\bold{u}_{x})
\end{bmatrix}$
& 
$\begin{bmatrix}
\bold{b}_{x}+\bold{s}_{x}^{k} \\
\bold{b}_{z}+\bold{s}_{z}^{k}
\end{bmatrix}-\begin{bmatrix}
\omega^2\bs{\rho}\bold{u}_{x}\\
\omega^2\bs{\rho}\bold{u}_{z}\end{bmatrix} $                                                                                                                                                                                     \\ \hline
$(\bold{V}_{P}^{2},\bold{V}_{S}^{2})$   
& 
$\begin{bmatrix}
\diag(\partial_{xx}\bold{u}_{x}\circ \bs{\rho})+\diag(\partial_{xz}\bold{u}_z \circ \bs{\rho}) &  \diag(\partial_{zz}\bold{u}_{x} \circ\bs{\rho})-\diag(\partial_{xz}\bold{u}_z \circ \bs{\rho}) \\
\diag(\partial_{zz}\bold{u}_z \circ\bs{\rho})+\diag(\partial_{xz}\bold{u}_{x}\circ\bs{\rho}) & \diag(\partial_{xx}\bold{u}_z \circ \bs{\rho})-\diag(\partial_{xz}\bold{u}_{x}\circ \bs{\rho})
\end{bmatrix} $ 
& 
$\begin{bmatrix}
\bold{b}_{x}+\bold{s}_{x}^{k} \\
\bold{b}_{z}+\bold{s}_{z}^{k}
\end{bmatrix}-\begin{bmatrix}
\omega^2\bs{\rho}\bold{u}_{x}\\
\omega^2\bs{\rho}\bold{u}_{z}\end{bmatrix} $
\\  \hline
\end{tabular}
\label{Table:L_lame_rho}
\end{table*}
%
%
%
%
%
%
%
%

%
%
%

\subsection{Dual ascent step}
To mitigate the instability issue during early iterations, a damped multiplier update can be employed, as suggested by \citet{gholami2023multiplier}. This update can be applied by replacing the multiplier update in equation \ref{ADMM_eta} with the following damped update:
\begin{equation}\label{eq:damped_multipliers}
\bold{s}^{k+1} = (\frac{k}{k+\xi})\left(\bold{s}^{k}+\bold{b}-\A(\bold{m}^{k+1})\bold{u}^{k+1}\right), 
\end{equation} 
where $\xi$ is a damping factor, typically chosen to be larger than 1. The damping factor can be adjusted based on the stability of the algorithm and the convergence behavior observed during the iterations. A larger damping factor can help stabilize the algorithm, but it may slow down convergence. On the other hand, a smaller damping factor may speed up convergence, but it can also introduce instability. 
By using the damped multiplier update, the ADMM algorithm can handle rough initial models more effectively and reduce the need for successive restarts, leading to improved convergence behavior \citep{gholami2023multiplier}.
\section{Numerical Example}
In this section, the performance of the proposed algorithm is assessed through a set of numerical examples. The quality of the inversion results is quantified using several metrics, including the ``source residual", ``data residual", and "model error".

The ``source residual" measures the discrepancy between the modeled sources computed using the updated model $\m^{k+1}$ and the desired sources $\b$. It is computed as the Euclidean norm of the difference between the modeled sources and the desired sources. 
The "data residual", on the other hand, measures the mismatch between the observed data $\d$ and the modeled data $\Pu^{k+1}$. It is computed as the  Euclidean norm of the difference between the observed data and the modeled data. 
Finally, the ``model error" quantifies the deviation between the estimated model $\m^{k+1}$ and the true model $\m^{*}$. It is computed as the Euclidean norm of the difference between the estimated model and the true model. 

These metrics provide quantitative measures of the accuracy of the inversion results. A smaller source residual and data residual indicate a better fit between the modeled and desired sources and data, respectively. Similarly, a smaller model error indicates a closer resemblance between the estimated model and the true model.
By evaluating these metrics, we can assess the effectiveness and reliability of the proposed algorithm in capturing the subsurface properties and recovering the true model. Moreover, the PDE operator, $\bold{A}$, is discretized using the optimal 9-point finite difference stencil proposed by \citet{chen2016modeling}. 

\subsection{Double circular toy model}
\subsubsection{On the role of parameterization}
\begin{figure}[!ht]
\centering 
\includegraphics[scale=.75]{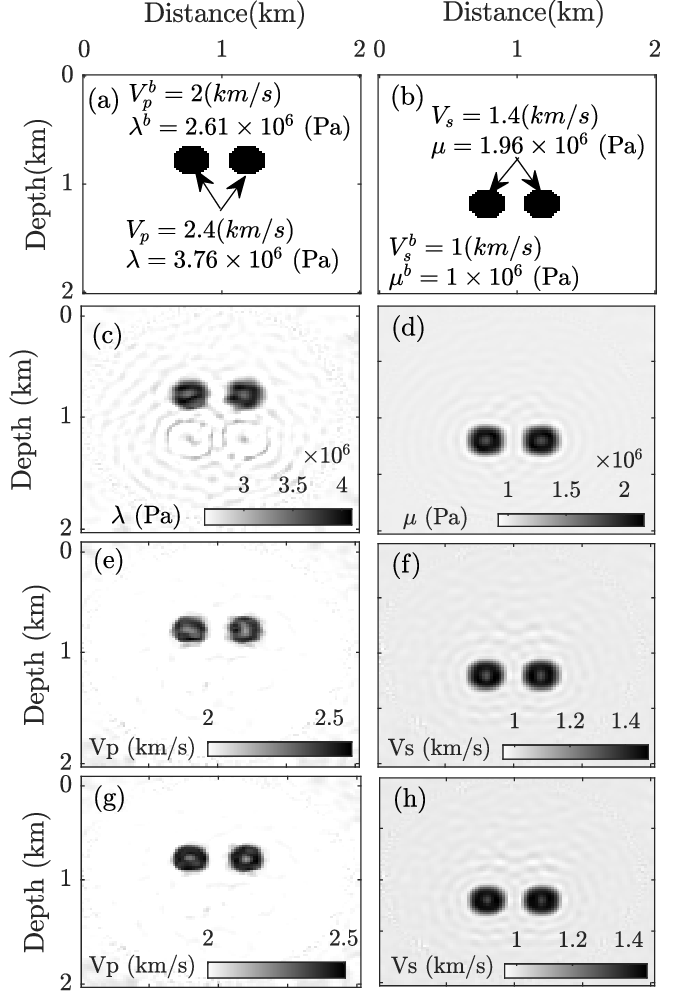}
\caption{Double circular model. (a-b) True ($\bold{V}_\text{P}$, $\bold{V}_\text{S}$) and ($\bs{\lambda}, \bs{\mu}$) models, in which parameters with superscripts (b) indicate the background values. The estimated models obtained with parameterization ($\bs{\lambda}, \bs{\mu}$) (c-d) , ($\bold{V}_\text{P}$, $\bold{V}_\text{S}$) (e-f) and ($\bold{V}_\text{P}^2, \bold{V}_\text{S}^2$) (g-h) after 70 iterations.}
\label{Fig:model_parameterization}
\end{figure}

In the assessment of the model parameterization for constant density elastic media, double circular heterogeneities are considered for both the ($\bold{V}_\text{P}$, $\bold{V}_\text{S}$) and ($\bs{\lambda}$, $\bs{\mu}$) parameterizations. These heterogeneities are located at different positions and embedded in a homogeneous background model, as shown in Figure~\ref{Fig:model_parameterization}a and Figure~\ref{Fig:model_parameterization}b, respectively.
It is important to note that for this analysis, velocity models were not constructed using Lam$\acute{\text{e}}$ parameters. The primary goal here is to investigate the potential cross-talk between the different parameterizations and assess their impact.
The circular acquisition setup (with radius of 1~km) consists of 16 vertically directed forces, denoted as $\bold{b}_{z}$, which emit a Ricker wavelet with a dominant frequency of 10 Hz. The wavefields are recorded by 128 two-component receivers.
To prevent any reflection artifacts, absorbing boundary conditions are implemented along the four sides of the modeling domain.
%
%
%
%
%
%
%
%
%
\begin{figure}[!h]
\centering 
\includegraphics[scale=.75]{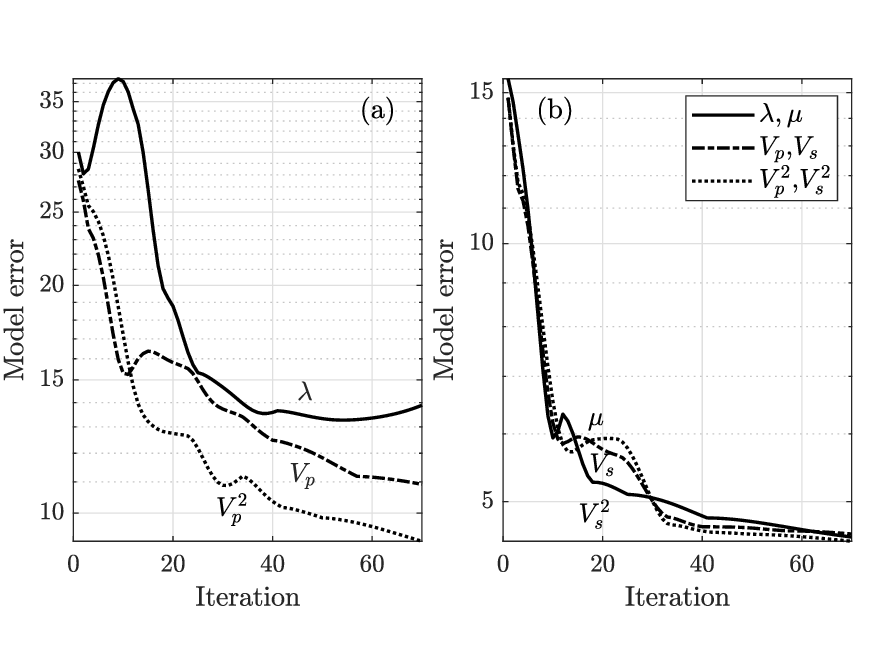}
\caption{Double circular model. The evolution of the computed model errors during iteration for different model parameterizations: (a) ($\bs{\lambda}, \bold{V}_\text{P}, \bold{V}_\text{P}^{2}$), and (b) ($\bs{\mu}, \bold{V}_\text{S}, \bold{V}_\text{S}^{2}$).}
\label{Fig:error_parameterization}
\end{figure}
%
%
%
%
%
%

\begin{figure}[!h]
\centering 
\includegraphics[scale=.65]{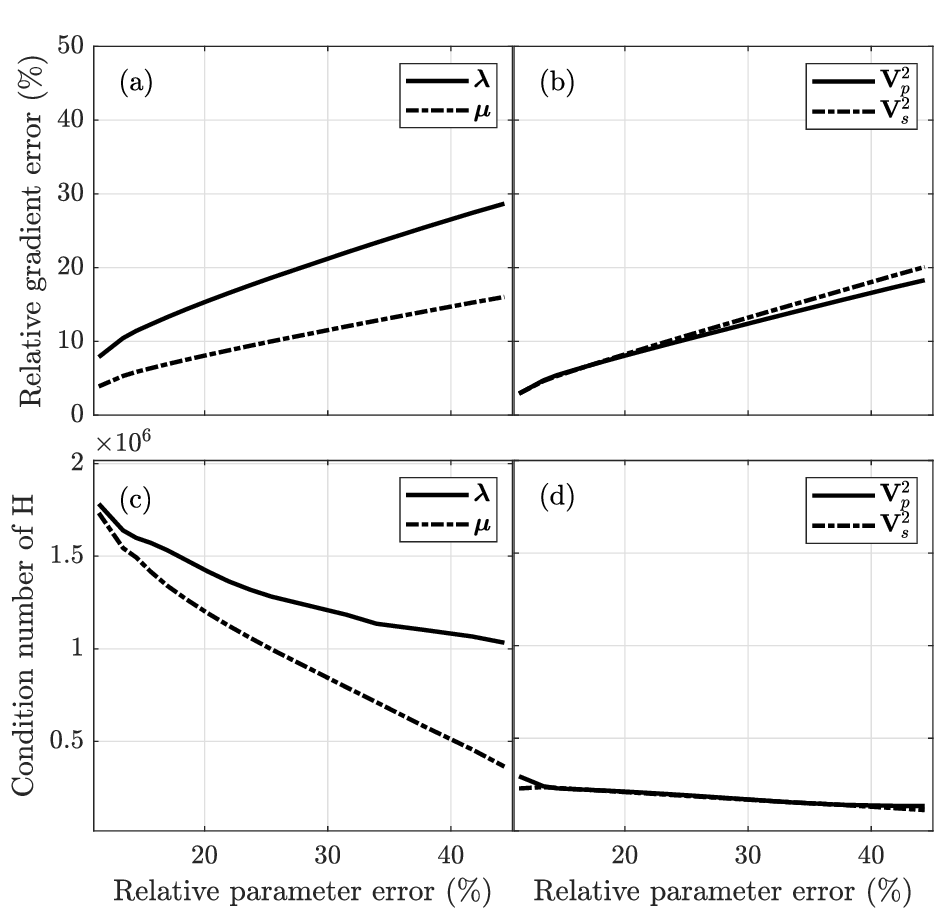}
\caption{Sensitivity analysis for $(\bs{\lambda}, \bs{\mu})$ and $(\bold{V}_\text{P}^2, \bold{V}_\text{S}^2)$ parameterizations. Relative error in the gradient vector versus the error in each parameter for (a) $(\bs{\lambda}, \bs{\mu})$ and (b) $(\bold{V}_\text{P}^2, \bold{V}_\text{S}^2)$. The condition number of the Hessian matrix versus relative error of the input parameters for (c) $(\bs{\lambda}, \bs{\mu})$ and (d) $(\bold{V}_\text{P}^2, \bold{V}_\text{S}^2)$.}
\label{Fig:error_rhs}
\end{figure}
In the analysis, background models were used to perform inversion for two frequencies, specifically 2.5 Hz and 5 Hz, concurrently. Three different parameterizations were employed: ($\bs{\lambda}$, $\bs{\mu}$), ($\bold{V}_\text{P}$, $\bold{V}_\text{S}$), and ($\bold{V}_\text{P}^2$, $\bold{V}_\text{S}^2$).
The results of the inversion after 70 iterations are shown in Figure~\ref{Fig:model_parameterization}c-d for the ($\bs{\lambda}$, $\bs{\mu}$) parameterization, Figure~\ref{Fig:model_parameterization}e-f for the ($\bold{V}_\text{P}$, $\bold{V}_\text{S}$) parameterization, and Figure~\ref{Fig:model_parameterization}g-h for the ($\bold{V}_\text{S}^2$, $\bold{V}_\text{S}^2$) parameterization. Generally, the inverted parameters obtained using ($\bold{V}_\text{P}$, $\bold{V}_\text{S}$) and ($\bold{V}_\text{P}^2$, $\bold{V}_\text{S}^2$) exhibit comparable accuracy.
During the early iterations of the ADMM algorithm, all parameterizations exhibit cross-talk, where the $\bold{V}_\text{S}$ parameter influences $\bold{V}_\text{P}$ and $\bs{\mu}$ influences $\bs{\lambda}$. However, as the iterations progress, these cross-talk effects are suppressed. There is a slight parameter leakage of $\bs{\mu}$ into $\bs{\lambda}$, but it can be regarded as noise and can be controlled through appropriate regularization techniques.
Figure~\ref{Fig:error_parameterization} illustrates the evolution of model errors versus iteration. It is worth noting that there is no significant parameter leakage observed in terms of velocity parameterization. This may be attributed to the influence of the Gauss-Newton Hessian matrix ($\bold{L}^T\bold{L}$), which will be discussed in the next section.
\subsubsection{Sensitivity Analysis}
In the analysis of model parameterization, the sensitivity of each parameter class is assessed to determine the more stable case between ($\bs{\lambda}$, $\bs{\mu}$) and ($\bold{V}_\text{P}^2$, $\bold{V}_\text{S}^2$). The sensitivity analysis involves solving a linear system associated with each parameter class and perturbing one variable by a specific amount to analyze the resulting relative error in the gradient vector $\bold{L}^T\bold{y}$ and the condition number of the corresponding Hessian matrix $\bold{L}^T\bold{L}$.
For this analysis, a surface acquisition geometry is considered. The $\bold{V}_\text{S}$ model is constructed using a constant Poisson's ratio of 0.24, while the Lamé parameters are constructed in the same precise manner.

Figure~\ref{Fig:error_rhs} displays the calculated gradient error and condition number for a frequency of 3 Hz. The results show that the gradient vector is more sensitive to errors in $\bs{\lambda}$ in the case of the ($\bs{\lambda}$, $\bs{\mu}$) parameterization. In contrast, the sensitivity is almost the same for the parameters in ($\bold{V}_{P}^2$, $\bold{V}_\text{S}^2$). Furthermore, the condition number of the Hessian matrix is an order of magnitude higher for the ($\bs{\lambda}$, $\bs{\mu}$) parameterization, indicating a poorly-conditioned parameterization.
Based on these observations, it is concluded that the ($\bold{V}_{S}^2$, $\bold{V}_{P}^2$) parameterization is more stable. 
This aligns with the observations made by \citet{kohn2012influence} and \citet{pan2018elastic}. 
Therefore, it is chosen as the more reliable parameterization for the subsequent examples.
%
%
\subsubsection{On the role of Hessian}
In order to investigate the role of the Gauss-Newton Hessian matrix ($\bold{L}^T\bold{L}$) in suppressing inter-parameter cross-talk, we consider the model update equation for the parameterization $(\bold{m}_p, \bold{m}_s) \equiv (\bold{V}_\text{P}^2, \bold{V}_\text{S}^2)$, which can be written as:
\begin{subequations}\label{eq:Hess_decomp}
\begin{align}
\bold{m}_{p} = \bold{H}_{11}^{-1}\bold{g}_{x}+\bold{H}_{12}^{-1}\bold{g}_{z}, \label{mp_update} \\
\bold{m}_{s} = \bold{H}_{21}^{-1}\bold{g}_{x}+\bold{H}_{22}^{-1}\bold{g}_{z}, \label{ms_update} 
\end{align}
\end{subequations}
where $\bold{H}^{-1}_{ij}$ represents the blocks of the inverse Hessian matrix and $\bold{g} = (\bold{g}_x, \bold{g}_z)$ is the gradient vector. 
It is crucial to note that the Hessian blocks $\bold{H}^{-1}_{ij}$ are diagonal and amenable to explicit computation.
%
\begin{figure}[!ht]
\centering 
\includegraphics[scale=.75,trim={0cm 1.5cm 0cm 0cm},clip]{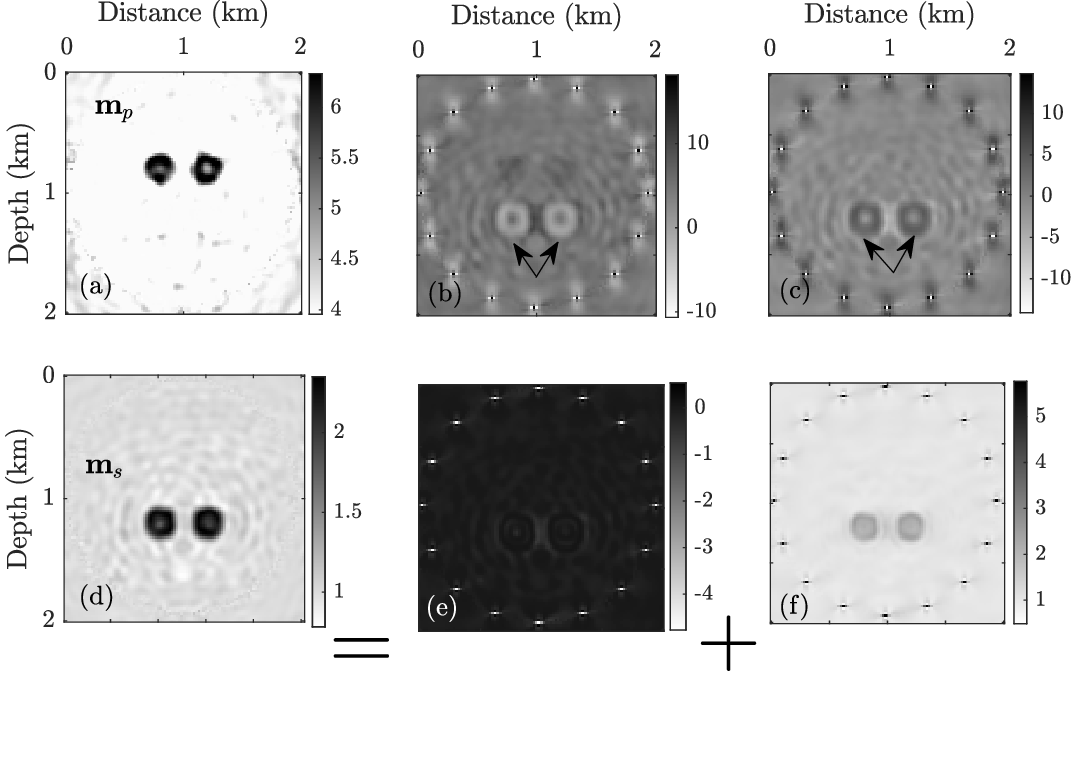}
\caption{On the role of the Hessian matrix. The model updates for $\bold{m}_{p}$  (a-c) and for $\bold{m}_{s}$ (d-f) at 20th iterations.
(a,d) are the total updates which can be decomposed into the parts due to each block of the Hessian; (b,f) diagonal blocks, and (c,e) off diagonal blocks. The diagonal blocks contain dominant features.
 }
\label{Fig:Hessian_role}
\end{figure}
Figure~\ref{Fig:Hessian_role} highlights the importance of the Hessian matrix in suppressing parameter cross-talk by demonstrating the influence of each Hessian block in the model update at the 20th iteration. From Figure~\ref{Fig:Hessian_role}b-c, it can be observed that there is considerable leakage of $\bold{m}_s$ into $\bold{m}_p$ (marked by black arrows). However, their destructive summation effectively suppresses this leakage. This toy experiment emphasizes the effectiveness of the Hessian in reducing parameter cross-talk during the inversion process. Additionally, the sparsity of the Hessian matrix makes it suitable for large-scale FWI problems.

\newpage
\subsection{SEG/EAGE overthrust model}
\begin{figure*}[h] 
\centering 
\includegraphics[scale=.75]{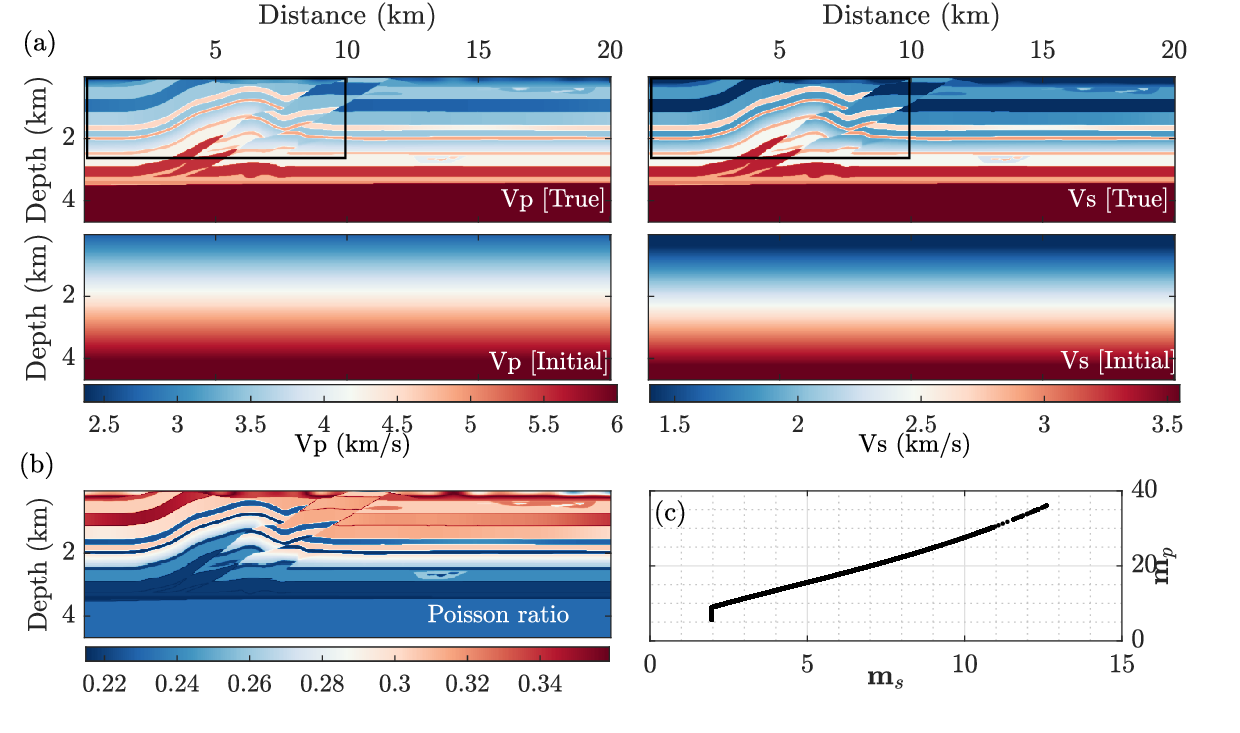}
\caption{(a) True and initial P-and S-wave velocity models, (b) Poisson's ratio, and (c) cross-plot of $\bold{m}_{p}$ versus $\bold{m}_{s}$.}
\label{Fig:over_initial}
\end{figure*}

\noindent
In the subsequent experiments, we focus on a 2D section of the 3D SEG/EAGE overthrust $\bold{V}_\text{P}$ model. The model dimensions are 4.67 km $\times$ 20 km with a grid interval of 25 m. We infer the $\bold{V}_\text{S}$ model from the $\text{V}_\text{P}$ model using an empirical relation proposed by \citet{brocher2005empirical}:
\begin{equation}
\begin{aligned}
\bold{V}_\text{S}=&0.7858-1.2344\bold{V}_\text{P}+0.7949\bold{V}_\text{P}^2 
-0.1238\bold{V}_\text{P}^3+0.0044\bold{V}_\text{P}^4.
\end{aligned}
\end{equation}
This relationship leads to higher values of Poisson's ratio near the surface, which requires dense spatial sampling in the forward modeling process. To avoid this requirement, we set the minimum value of the $\bold{V}_\text{S}$ model to be 1.4 km/s (Figure~\ref{Fig:over_initial}a).
However, the inversion task remains challenging due to the presence of multiple Poisson's ratios within the medium. In the context of EFWI, Poisson's ratio plays a crucial role in influencing the accuracy and stability of the inversion process \citep{xu20142d, brossier2009seismic}.

The surface acquisition setup for the inversion process involves 134 sources spaced at intervals of 150 m. Ricker wavelets with dominant frequencies of 10 Hz are used as the source wavelets ($\bold{b}_x$ and $\bold{b}_z$). The receivers consist of 401 two-component sensors positioned every 50 m. Absorbing boundary conditions are used on top of the model, unless otherwise stated.
The inversion stage begins with models that linearly increase with depth. The $\bold{V}_\text{S}$ model ranges from 1.2 km/s to 3.8 km/s, while the $\bold{V}_\text{P}$ model ranges from 2.7 km/s to 6.5 km/s (Figure~\ref{Fig:over_initial}a).
The computed Poisson's ratio and the cross-plot of $\bold{m}_s-\bold{m}_p$ are shown in Figures~\ref{Fig:over_initial}b and c, respectively.
The inversion is carried out in three cycles following the standard multiscale strategy \citep{bunks1995multiscale}, with individual frequencies ranging from 3 Hz to 13 Hz in intervals of 0.5 Hz. The cycles are as follows:

\begin{enumerate}[topsep=1pt, partopsep=1pt]
\item Cycle 1: Frequencies range from  3~Hz to 6~Hz.
\item Cycle 2: Frequencies range from  3~Hz to 7.5~Hz.
\item Cycle 3: Frequencies range from  3~Hz to 13~Hz.
\end{enumerate}

The inversion process consists of 20 iterations for the first frequency component and 10 iterations for the subsequent frequencies, resulting in a total of 410 iterations. A fixed penalty parameter of $\beta = 10^6$ is used throughout the inversion process.
%
%
%
%
%
%
\begin{figure}[!h] 
\centering 
\includegraphics[scale=.7]{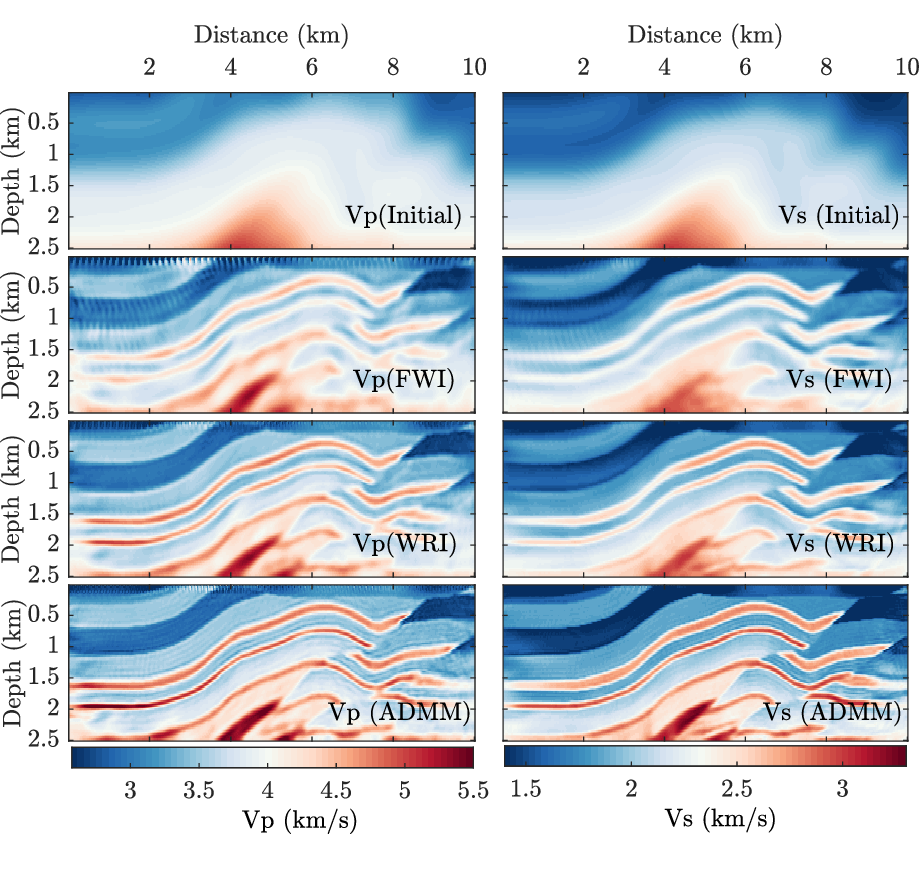}
\caption{Comparison between classic FWI, WRI and ADMM for the case of smooth initial
models.}
\label{Fig:over_FWI_vs_ADMM_smooth}
\end{figure}
%
%
\begin{figure}[!h] 
\centering 
\includegraphics[scale=.7]{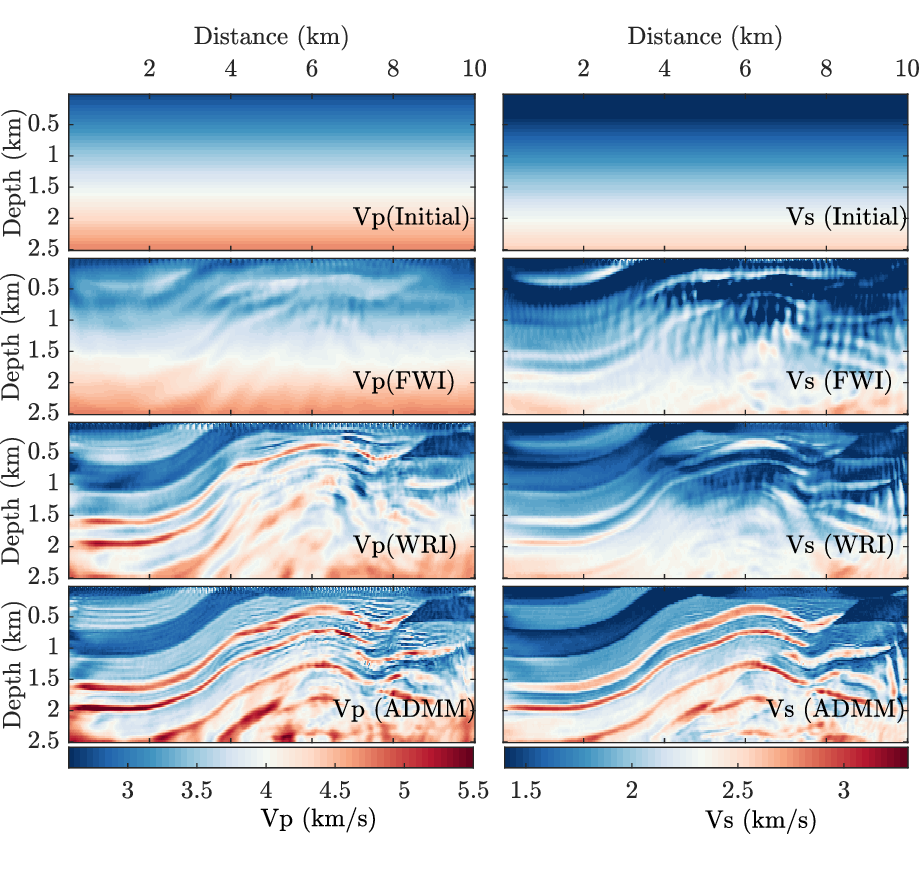}
\caption{Comparison between classic FWI, WRI and ADMM for the case of linearly
increased initial model.}
\label{Fig:over_FWI_vs_ADMM_grad}
\end{figure}
%
%
\begin{figure}[!h] 
\centering 
\includegraphics[scale=.7]{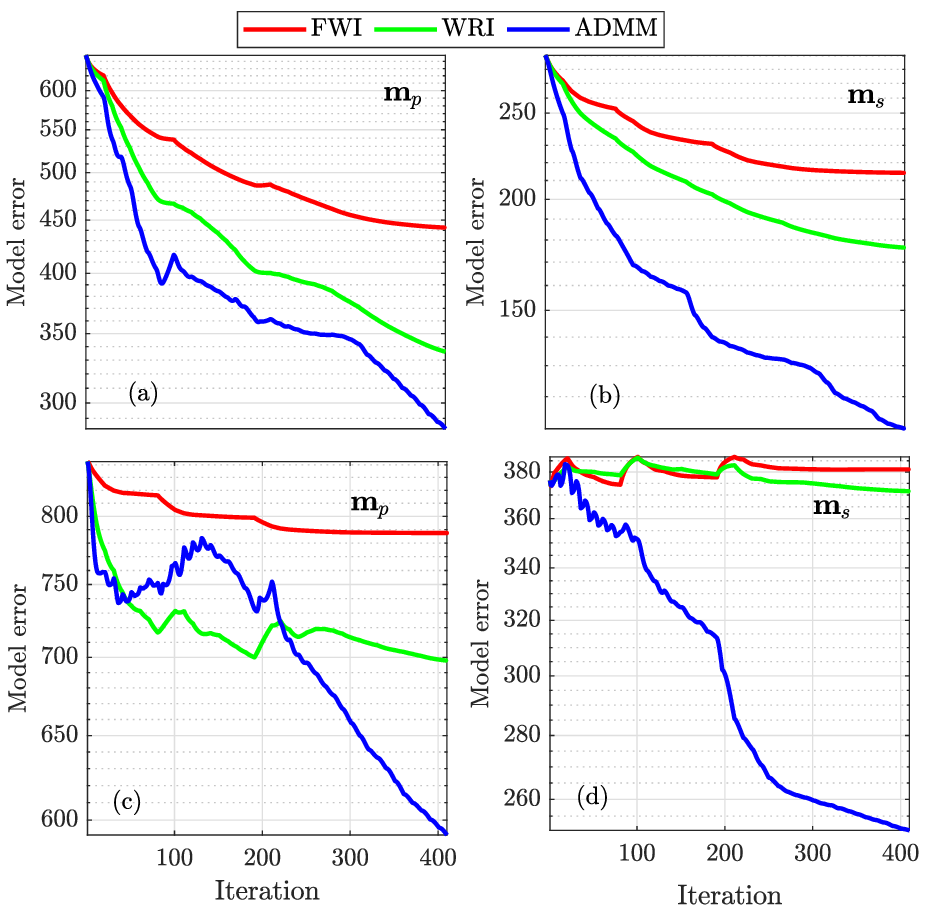}
\caption{Evolution of the model error versus iteration for classic FWI, WRI and ADMM started with (a,b) smooth initial models (Figure~\ref{Fig:over_FWI_vs_ADMM_smooth}) and (c,d) linearly increased initial models (Figure~\ref{Fig:over_FWI_vs_ADMM_grad}).}
\label{Fig:over_FWI_vs_ADMM_error}
\end{figure}
%
\begin{figure}[!h] 
\centering 
\includegraphics[scale=.7,trim=0cm 2.5cm 0cm 0cm,clip]{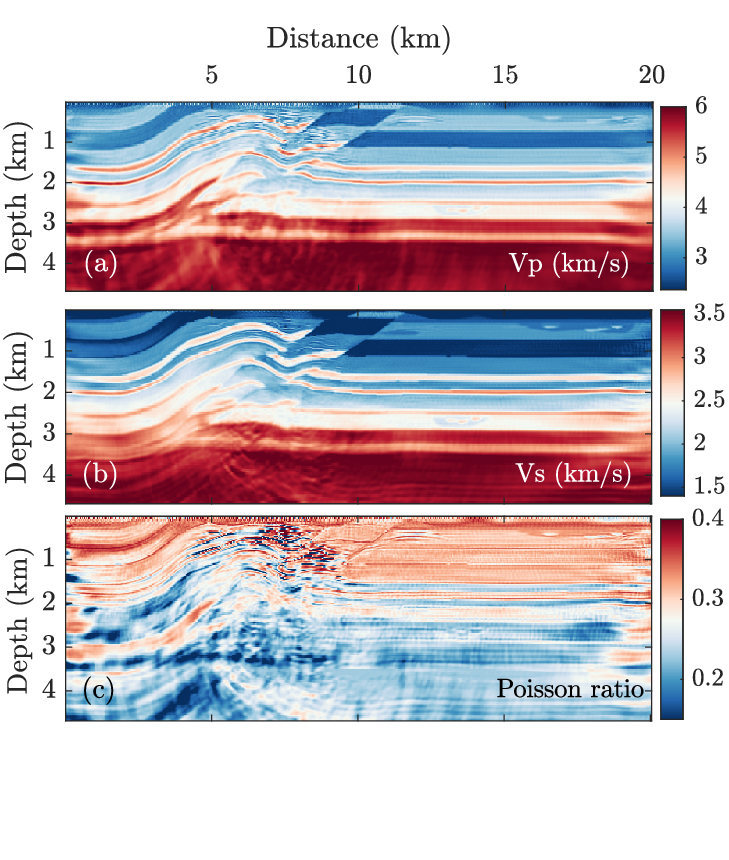}
\caption{Estimated (a) Vp, (b) Vs, and (c) Poisson ratio using squared velocity parametrization. }
\label{Fig:over_ALF}
\end{figure}
First, we conducted a comparative analysis to assess the efficacy of the proposed method against traditional FWI. The evaluation focused on a specific segment of the SEG/EAGE overthrust model, demarcated by a rectangle in Figure~\ref{Fig:over_initial}a. Two distinct starting models were utilized: one derived by smoothing the original models (shown at the top of Figure \ref{Fig:over_FWI_vs_ADMM_smooth}), and the other featuring linear depth-dependent increments, as illustrated at the top of Figure~\ref{Fig:over_FWI_vs_ADMM_grad}. \\
Three inversion methods were employed: classical FWI based on the reduced formulation, the proposed algorithm leveraging ADMM, and WRI, where WRI was applied with Lagrange multipliers set to zero. Figure~\ref{Fig:over_FWI_vs_ADMM_smooth} shows the inversion results from the smooth initial models, with representations of initial models and reconstructed models for classical FWI, WRI, and ADMM displayed from top to bottom, respectively. Given the accuracy of the initial models, all three methods successfully converged to the true solution. Notably, the use of WRI and ADMM techniques led to enhanced resolution in the inverted models, with the ADMM exhibiting superior accuracy.
However, a substantial shift in the quality of estimated models occurred when using a less accurate starting model (Figure~\ref{Fig:over_FWI_vs_ADMM_grad}). Classic FWI became trapped in a local minimum, failing to reconstruct both $\bold{V}_\text{P}$ and $\bold{V}_\text{S}$ models. Conversely, the WRI approach improved reconstruction, although the $\bold{V}_\text{S}$ model still exhibited artifacts. Remarkably, the ADMM methodology outperformed other approaches, demonstrating superior performance.
Quantitative validation of these observations is presented in Figure~\ref{Fig:over_FWI_vs_ADMM_error}, where subfigures a and b illustrate the model error curves for $\bold{V}_\text{P}$ and $\bold{V}_\text{S}$ models with smoothed initial models. Subfigures c and d depict similar curves for cases with linearly increased initial models. These model error curves underscore the significant convergence rate of the ADMM technique for both scenarios.
\newpage
\newpage
\subsubsection{Elastic FWI with ADMM}
\begin{figure}[!h] 
\centering 
\includegraphics[scale=.7]{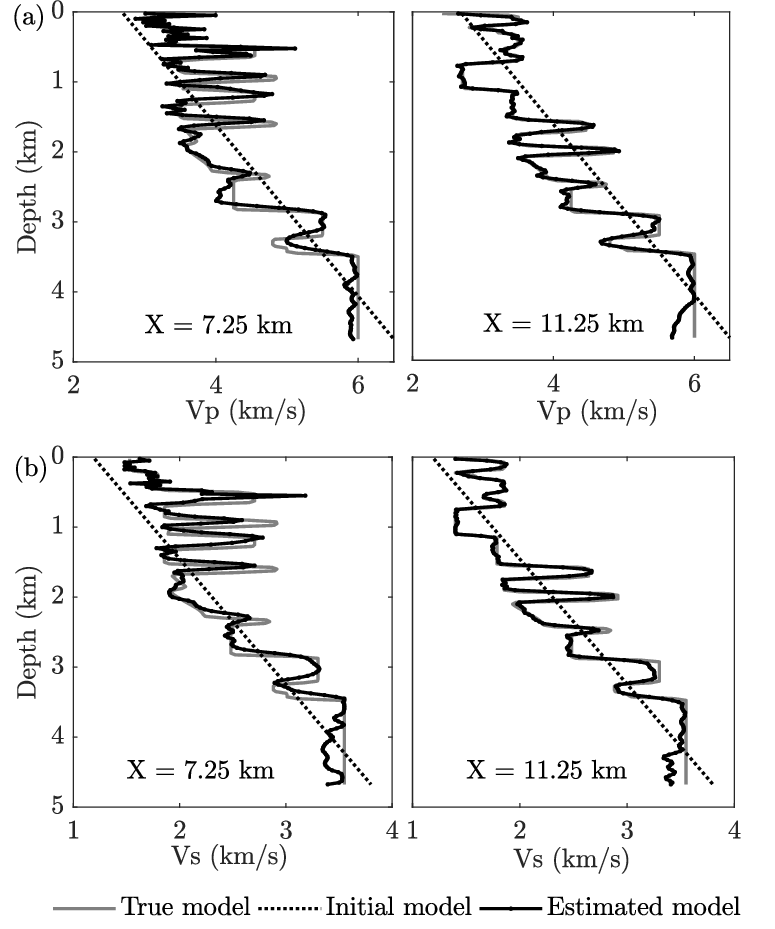}
\caption{Comparison between vertical velocity profiles at different locations labeled by $X$ for (a) $\bold{V}_\text{P}$ model and  (b) $\bold{V}_\text{S}$ model.}
\label{Fig:over_ALF_vs_vplog}
\end{figure}

We continued the investigation for EFWI using the ADMM method and investigated several characteristics to increase its performance. The inversion was performed to reconstruct models shown in Figure~\ref{Fig:over_initial}a.  Figures~\ref{Fig:over_ALF}a and b show the final inversion results for the Vp and Vs models, respectively. Figure~\ref{Fig:over_ALF}c displays the computed Poisson's ratio for the final estimated models.
To evaluate the performance of the method, Figure~\ref{Fig:over_ALF_vs_vplog} compares the extracted velocity profiles at various locations. 
%
%
%
%
%
%
%
%
%
%
%
%
%
%
%
%
%
\begin{figure}[!htb] 
\centering 
\includegraphics[scale=.6]{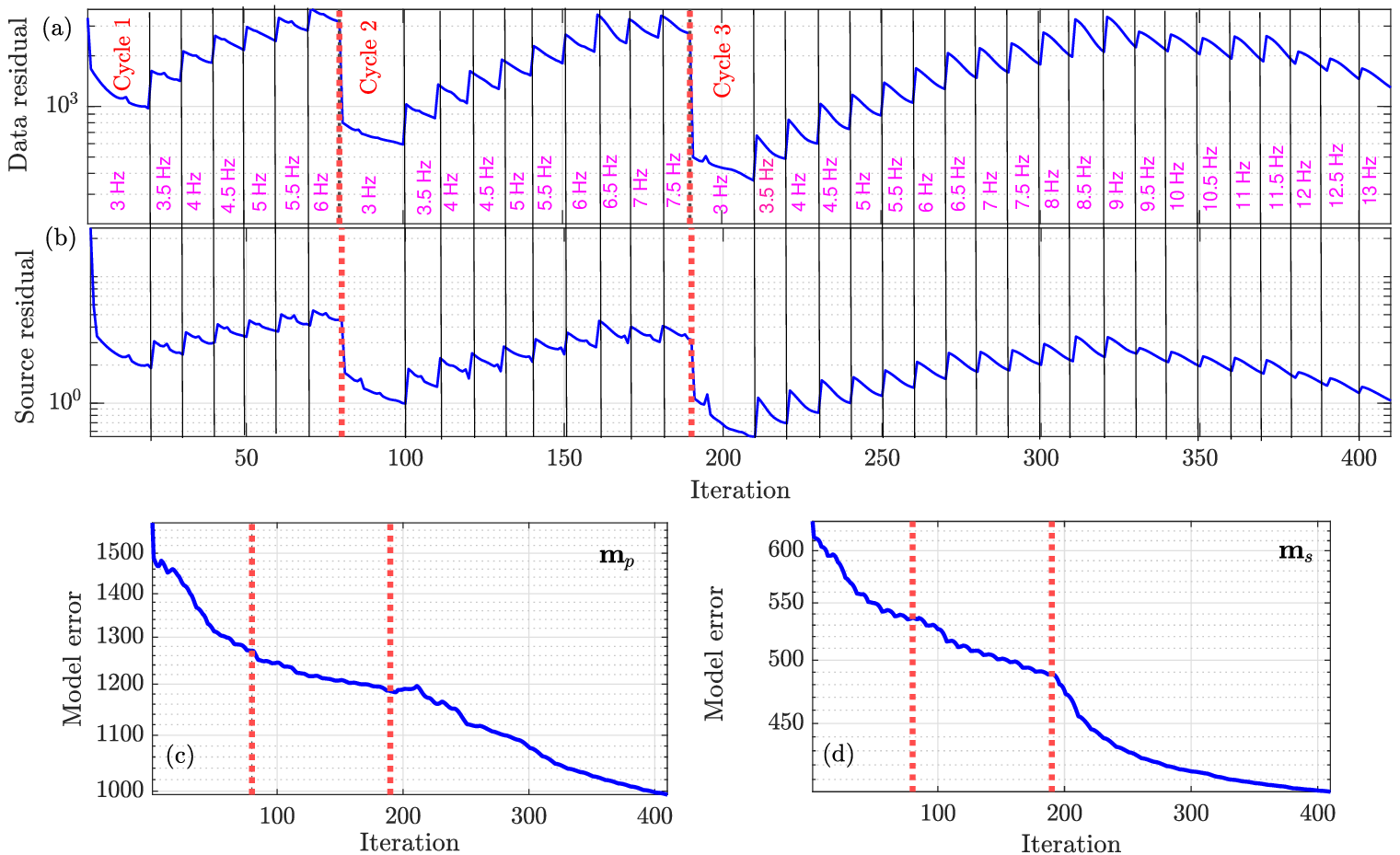}
\caption{Evolution of the (a) data residual, (b) source residual, (c) $\bold{m}_p$ model error, and (d) $\bold{m}_s$ model error versus inversion iteration.}
\label{Fig:over_ALF_vs_errors}
\end{figure}
%
%
\begin{figure}[!htb] 
\centering 
\includegraphics[scale=.5,trim=0cm 0.5cm 2cm 0cm,clip]{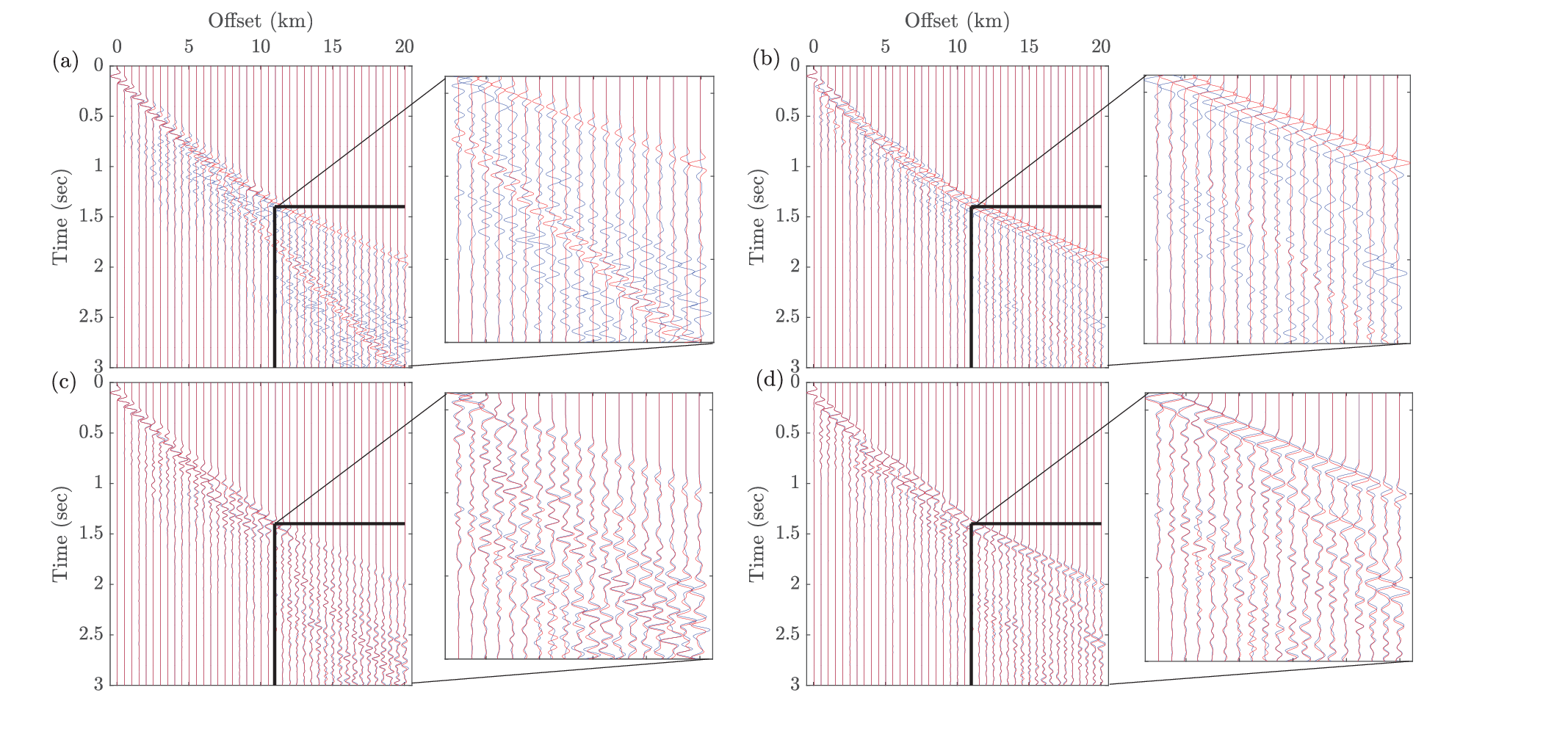}
\caption{Direct comparison of synthetic seismograms is shown in (a) and (c) for the horizontal component, and in (b) and (d) for the vertical component. Panels (a) and (b) show the seismograms computed from the true model ({blue}), overlays in a red wiggle plot representing the seismograms computed from the initial model. Panels (c) and (d) show the seismograms computed from the inverted model, also with overlays.}
\label{Fig:Shot_diff_1}
\end{figure}
The suggested method demonstrates good accuracy in reconstructing the models.
Figure~\ref{Fig:over_ALF_vs_errors} illustrates the convergence curves associated with the inversion process over the iterations. The figure clearly illustrates the decrease in errors for each frequency component throughout the iterations. It is noteworthy that, due to the distinct nature of the problems solved for each frequency, there might be instances of error increase during the transition from one frequency to another in the inversion process.
In Figure~\ref{Fig:Shot_diff_1} synthetic seismograms (horizontal and vertical components) are computed in the initial, true and estimated models to evaluate the fit of the phase and amplitude data and to determine the accuracy of the velocity reconstruction. To determine if the initial and true models were affected by cycle skipping, we superimpose the seismograms computed in the initial and estimated models ({in red}) over those computed in the true model ({in blue}). It can be seen that seismograms computed in the initial models are cycle skipped. However, a very good match between sesimograms computed in true and estimated models are obtained.
In the subsequent sections, we will explore the influence of various factors and assess the effectiveness of the algorithm.
\newpage
\subsubsection{The role of damping Lagrange multipliers}
\begin{figure*}[!h] 
\centering 
\includegraphics[scale=.65]{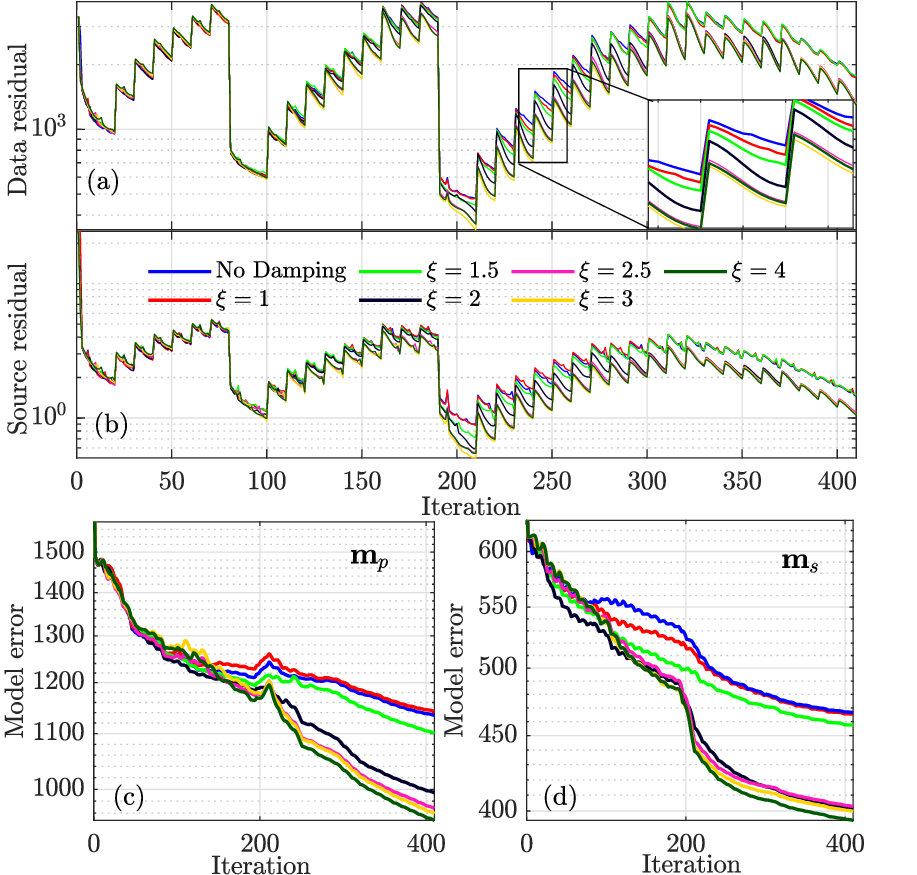}
\caption{The effect of the damping multipliers on the stability of the ADMM iterations (equation~\ref{eq:damped_multipliers}) using different values of $\xi$.  Evolution over iteration of data residual (a), source residual (b), $\bold{m}_p$ model error (c) and $\bold{m}_s$ model error (d).} 
\label{Fig:multipliers_damping}
\end{figure*}

\noindent
In order to address the issue of slow convergence caused by an inaccurate initial model, we applied damping to the Lagrange multiplier update in the early iterations of the ADMM algorithm. We tested various damping parameters ($\xi$) of 1, 1.5, 2, 2.5, 3, and 4, as described in equation \ref{eq:damped_multipliers}.
The performance of the ADMM iterations was evaluated based on the computed source residual, data residual, and model error over iterations. The quantitative analysis results are presented in Figure~\ref{Fig:multipliers_damping}.
The findings indicate that damping has a significant positive impact on the convergence of the algorithm, particularly for damping values of $\xi=2.5$, $3$, and $4$. After considering these results, we selected a damping value of $\xi=4$ for the subsequent experiments.

%
%
\subsubsection{The role of penalty parameter}
\begin{figure*}[!h] 
\centering 
\includegraphics[scale=.65]{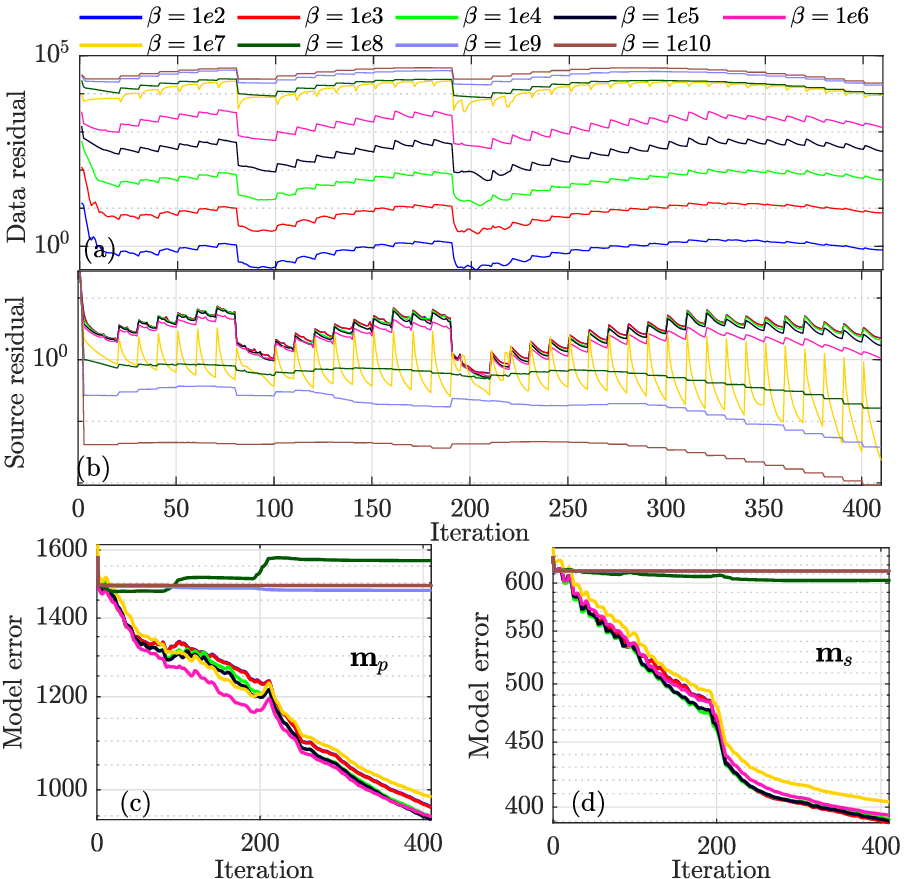}
\caption{On the choice of the penalty parameter: The performance of the AL based EFWI with different (fixed) values of the penalty parameter ($\beta$).  Evolution over iteration of data residual (a), source residual (b), $\bold{m}_p$ model error (c) and $\bold{m}_s$ model error (d).} 
\label{Fig:test_lambda}
\end{figure*}

\noindent
To assess the performance of the algorithm with different penalty parameter values, we conducted a test using values ranging from $\beta=10^2$ to $\beta=10^{10}$. We compared the evolution of the data residual, source residual, and model error for different values of $\beta$, as shown in Figure~\ref{Fig:test_lambda}.
The results indicate that when the value of $\beta$ is too low ($\beta<10^4$), the algorithm places more emphasis on minimizing the data term $\|\Pu - \d\|_2^{2}$. On the other hand, for large values of $\beta$ ($\beta>10^7$), the source term is given more weight. In both cases, the algorithm fails to produce satisfactory results.
However, penalty parameter values in the range of $\beta \in (10^4, 10^7)$ yielded satisfactory performance, balancing the influence of both the data and source terms. \\

%
\newpage
\subsubsection{Inversion of noise-contaminated data}
\begin{figure}[!htb]
\centering
\includegraphics[scale=.65,trim=0cm 0cm 0cm 0cm,clip]{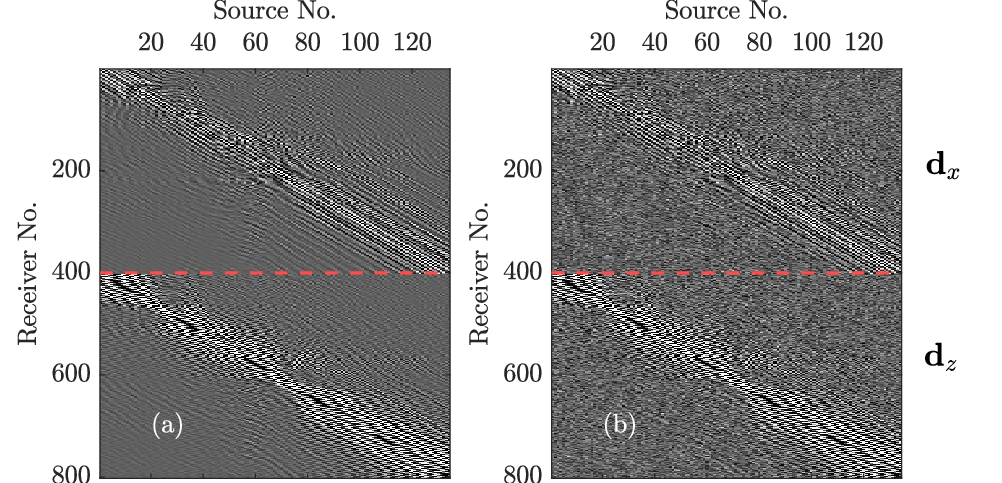}
\caption{Comparison between (a) noise-free and  (b) noise contaminated 7.5~Hz data with S/R of 10~dB. The dashed red line separates horizontal and vertical components of the receivers.}
\label{fig:over_data_noise}
\end{figure}

\begin{figure*}[!htb] 
\centering 
\includegraphics[scale=.6]{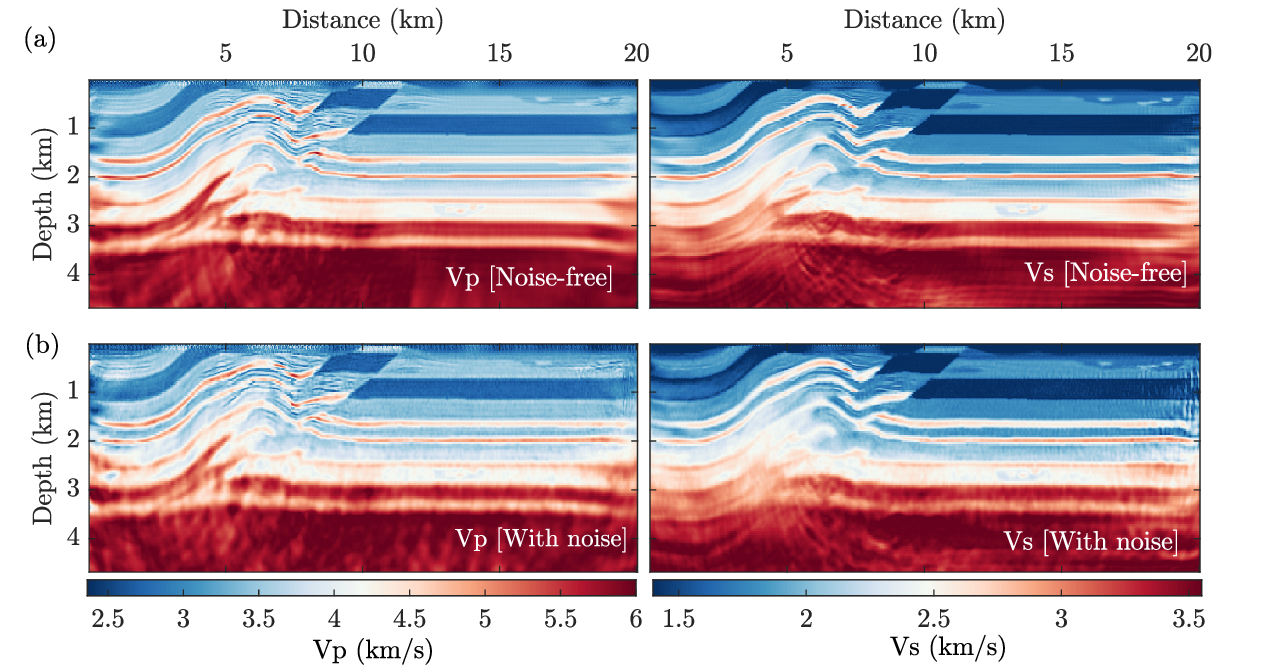}
\caption{Inversion results obtained from (a) noise-free data and (b) noisy data.}
\label{fig:over_inv_res_noise}
\end{figure*}

\noindent
In order to assess the effectiveness of our method in the presence of noise, we conducted an experiment using noisy data. Gaussian distributed random noise with a signal-to-noise ratio (S/R) of 10 dB was added to the data.
For a frequency component of 7.5 Hz, Figure~\ref{fig:over_data_noise} provides a comparison between the noise-free data and the noisy data in the source-receiver coordinate. Despite the added noise, the main features of the data are still preserved.
We then performed the inversion using the same setup, but with the noisy data. The resulting Vp and Vs models are displayed in Figure~\ref{fig:over_inv_res_noise}b. It can be observed that the impact of noise on the inversion results is not significant when compared to the results obtained from the noise-free data (Figure~\ref{fig:over_inv_res_noise}a).
This suggests that our method is relatively robust to noise and can produce reliable inversion results even in the presence of noise.
%
%
%
%
%
%
%
%
%
\subsubsection{The role of multi-physical constraints}
\begin{figure}[!htb]
\centering
\includegraphics[scale=.5]{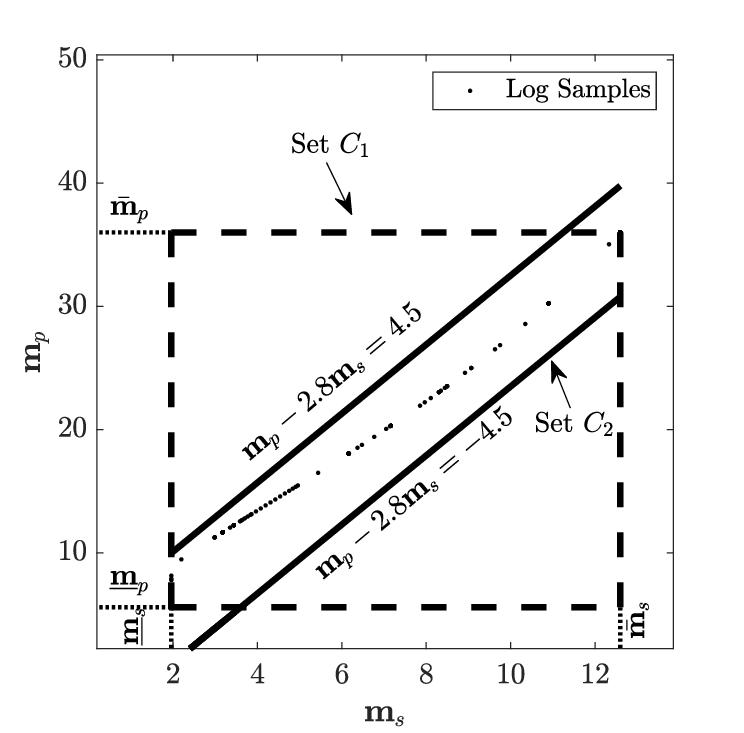}
\caption{Elastic FWI with physical constraints: ($\m_{s}-\m_{p}$) cross-plot (the black-dots) for a true vertical profile and two closed-convex sets ($C_1$ and $C_2$) defined as the physical constraints according to vertical velocity log samples.}
\label{Fig:Over_covex_sets}
\end{figure}

\begin{figure*}[!htb]
\centering
\includegraphics[scale=.65]{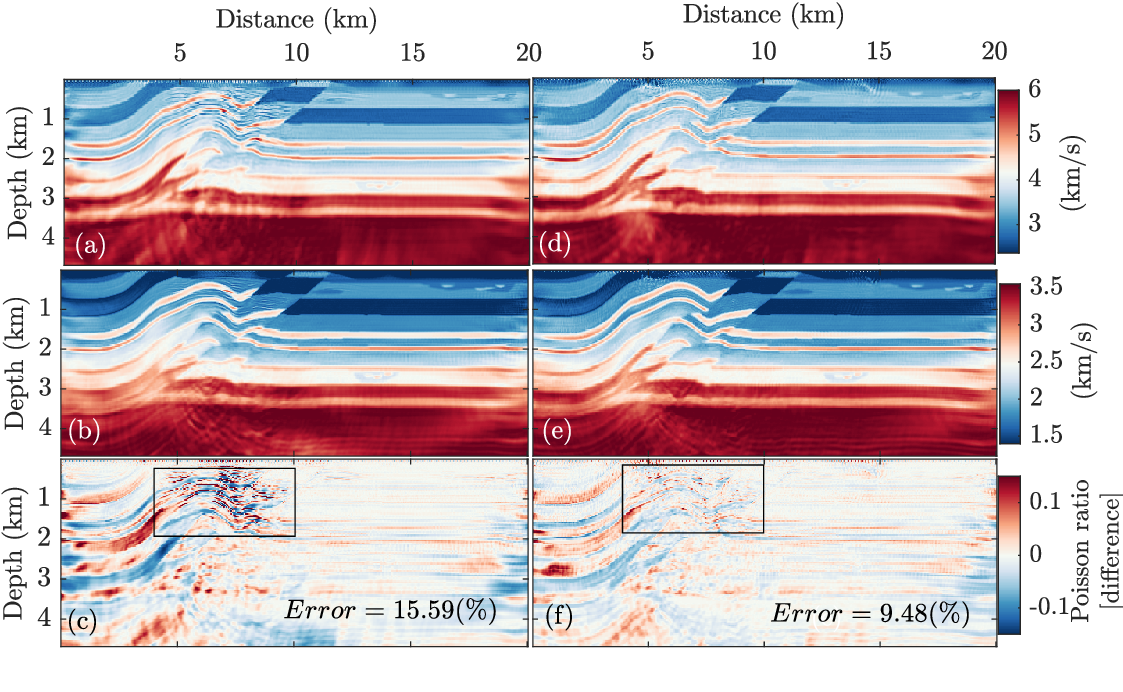}
\caption{Inversion with physical constraints. The (a) estimated Vp, (b) estimated Vs, and (c) difference between true and estimated Poisson's ratio when forcing $\m \in C_1$. (d-f) same as (a-c) for the case forcing $\m \in C_1 \cap C_2$.}
\label{fig:proj_inv_res}
\end{figure*}

\begin{figure*}[!htb]
\centering
\includegraphics[scale=.65,trim=0cm 0cm 1cm 0cm,clip]{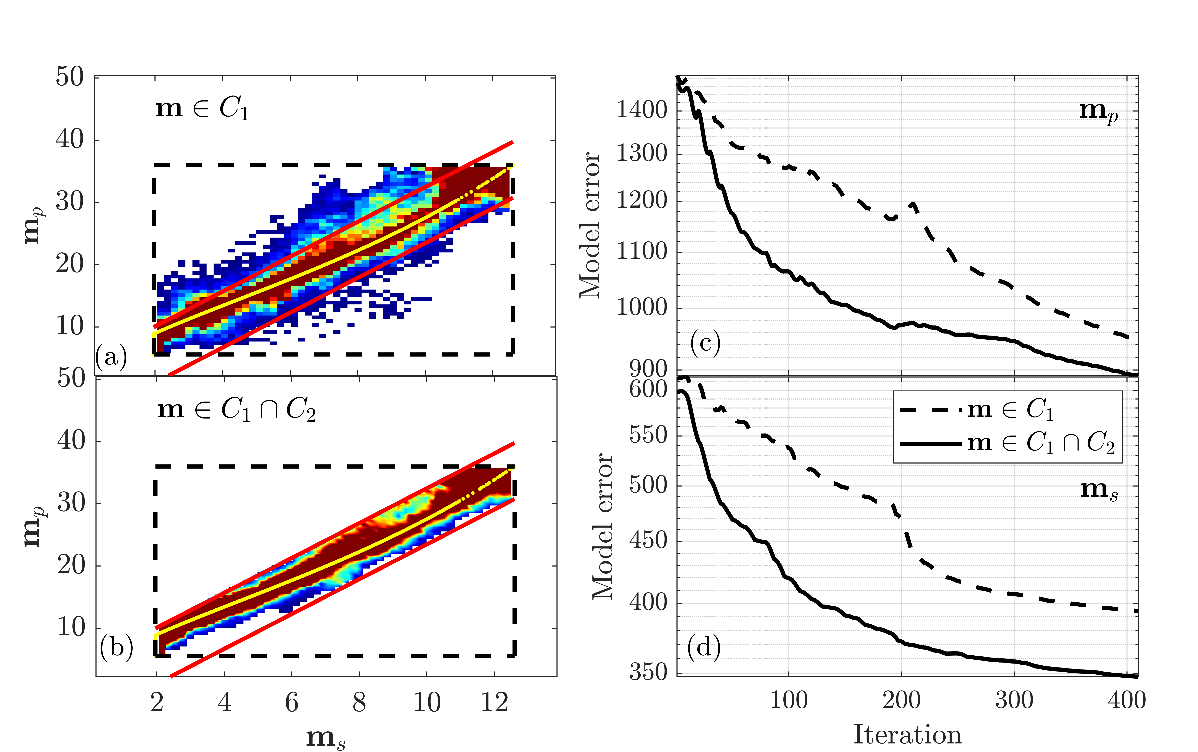}
\caption{Inversion with physical constraints. (a) Cross plot of $\m_{p}$ (Figure \ref{fig:proj_inv_res}a) versus $\m_{s}$ (Figure \ref{fig:proj_inv_res}b) obtained by forcing $\m \in C_1$.  (b) Cross plot of $\m_{p}$ (Figure \ref{fig:proj_inv_res}d) versus $\m_{s}$ (Figure \ref{fig:proj_inv_res}e) obtained by forcing $\m \in C_1 \cap C_2$, {where the yellow curve indicates the true cross plot.} The evolution of model errors over iteration for (c) $\m_{p}$ and (d) $\m_{s}$.}
\label{fig:proj_errors}
\end{figure*}

\noindent
In this subsection, we conducted experiments to test the impact of implementing physical constraints on the inversion results. We constructed two sets, $C_1$ and $C_2$, based on the true models (Figure~\ref{Fig:Over_covex_sets}).
$C_1$ was constructed based on the lowest and highest values of the true models (represented by dashed lines), while $C_2$ was constructed based on the cross-plot of Vp and Vs logs at $X=2.5$ km (represented by black dots).
We tested two cases: projection onto $C_1$ and projection onto $C_1 \cap C_2$. The inverted results for these cases are shown in Figure~\ref{fig:proj_inv_res}a, b, d, and e. To assess the quality of the estimates, we also plotted the Poisson's ratio residual sections.
It can be observed that when we implemented the constraint $\bold{m} \in C_1 \cap C_2$, the quality of the results significantly improved (Figure~\ref{fig:proj_inv_res}f versus c), especially in the faulted zone indicated by the black rectangle. This improvement is further supported by the cross-plots shown in Figure~\ref{fig:proj_errors}a, b, as well as the evolution of the model errors over the iterations in Figure~\ref{fig:proj_errors}c, d.
These results demonstrate that incorporating physical constraints into the inversion process, such as constraining the model to lie within the intersection of multiple convex sets, can greatly enhance the accuracy and reliability of the inversion results.

%
%
%
%
%
%
%
\subsubsection{Inversion with source sketching}
\noindent
In order to address the computational burden associated with solving the PDE in the ADMM iteration, we propose the use of source sketching, a method based on randomized discrete cosine transform (DCT) as the sketching matrix. The theoretical background of this method can be found in \citep{Aghazade_2021_RSS}.
%
\begin{figure}[!h]
\centering
\includegraphics[scale=.55,trim=0cm 0cm 0cm 0cm,clip]{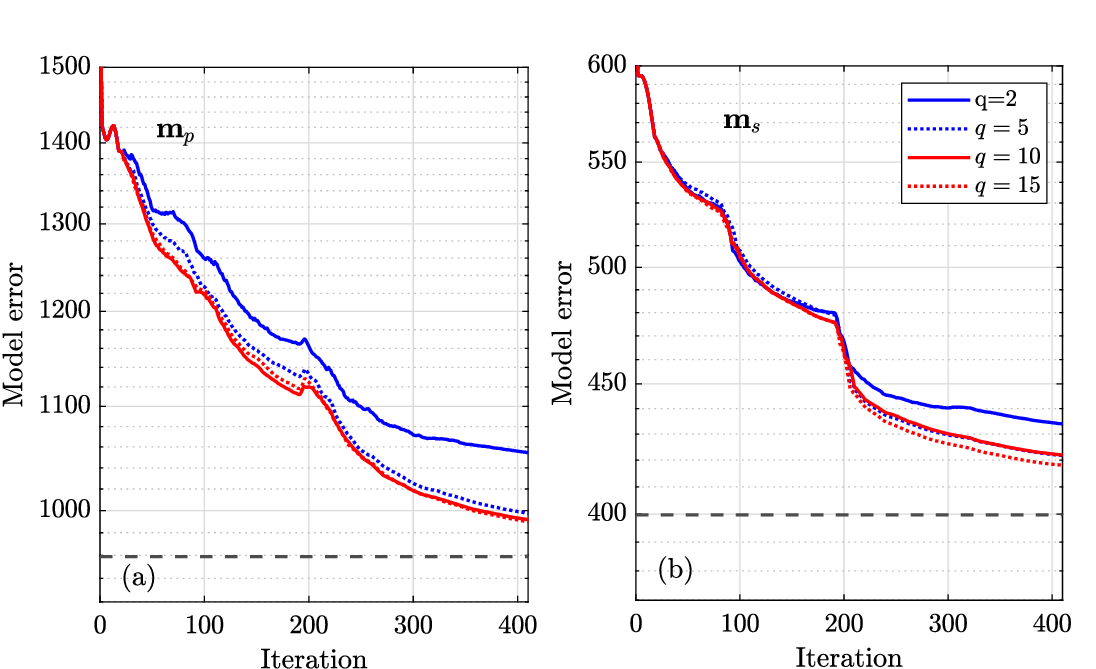}
\caption{Elastic FWI with source sketching: The evolution of the model error during iteration using different sketched sources for (a) $\bold{m}_p$ and  (b) $\bold{m}_s$.  The inversion is initialized with 50 sketched sources up to 20 iterations followed by sketched sources of $q=2, 5, 10, 15$. The horizontal dashed line indicates the final model error obtained by the full source inversion (134 sources).}
\label{fig:test_sketching}
\end{figure}

To assess the performance of the source sketching method, we conducted a series of experiments with different numbers of sketched sources. The total number of physical sources in the inversion was 134, and a total of 410 iterations were performed. During the initial 20 iterations, we used 50 sketched sources instead of the full set of 134 sources. For the remaining 390 iterations, we conducted experiments with different numbers of sketched sources: $q$=2, 5, 10, and 15. 
The evolution of the model error for each experiment is shown in Figure~\ref{fig:test_sketching}. The horizontal dashed line represents the final model error achieved through full-source (134 sources) inversion. In terms of accuracy, we observed that as the number of sketched sources increased, the final inversion results approached those obtained with the full set of sources. Regarding computational efficiency, the number of PDEs solved was significantly reduced compared with full-source inversion: from 54940 for the full set of sources to 1780 ($q=2$), 2950 ($q=5$), 4900 ($q=10$), and 6850 ($q=15$) for the sketched sources. Excluding the first 20 iterations, the computation speed-up, measured as $(1-q/134)\times 100$, was approximately 98.5\% ($q=2$), 96.2\% ($q=5$), 92.5\% ($q=10$), and 88.8\% ($q=15$), indicating the efficiency of the source sketching method.
\subsubsection{Inversion with free surface effects} \label{free_surface}
\begin{figure}[!h] 
\centering 
\includegraphics[scale=.65,trim=0cm 2.5cm 0cm 0cm,clip]{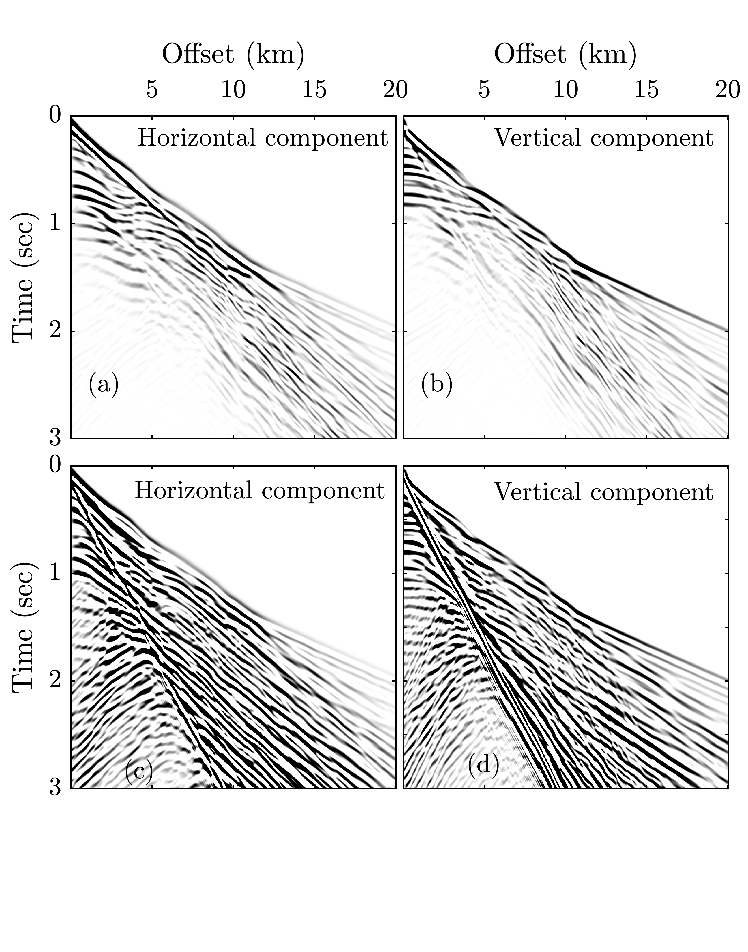}
\caption{Comparison between true seismograms computed  (a-b) without and (c-d) with free surface effects.}
\label{Fig:shot_FS}
\end{figure}

\begin{figure}[!h] 
\centering 
\includegraphics[scale=.6]{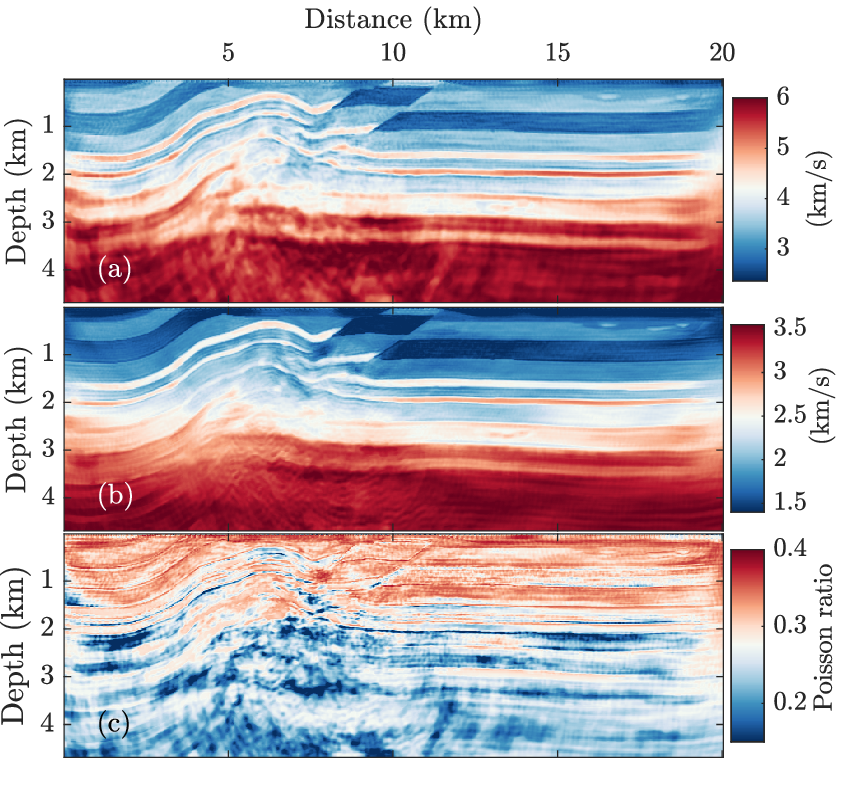}
\caption{Inversion with free surface effects. The (a) estimated Vp model, (b) estimated Vs model, and (c) {Estimated Poisson's ratio}. The inversion was performed by forcing $\bold{m}\in C_1\cap C_2$ and using the 1D initial models in Figure \ref{Fig:over_initial}.}
\label{Fig:inv_res_FS}
\end{figure}

\begin{figure}[!h] 
\centering 
\includegraphics[scale=.65]{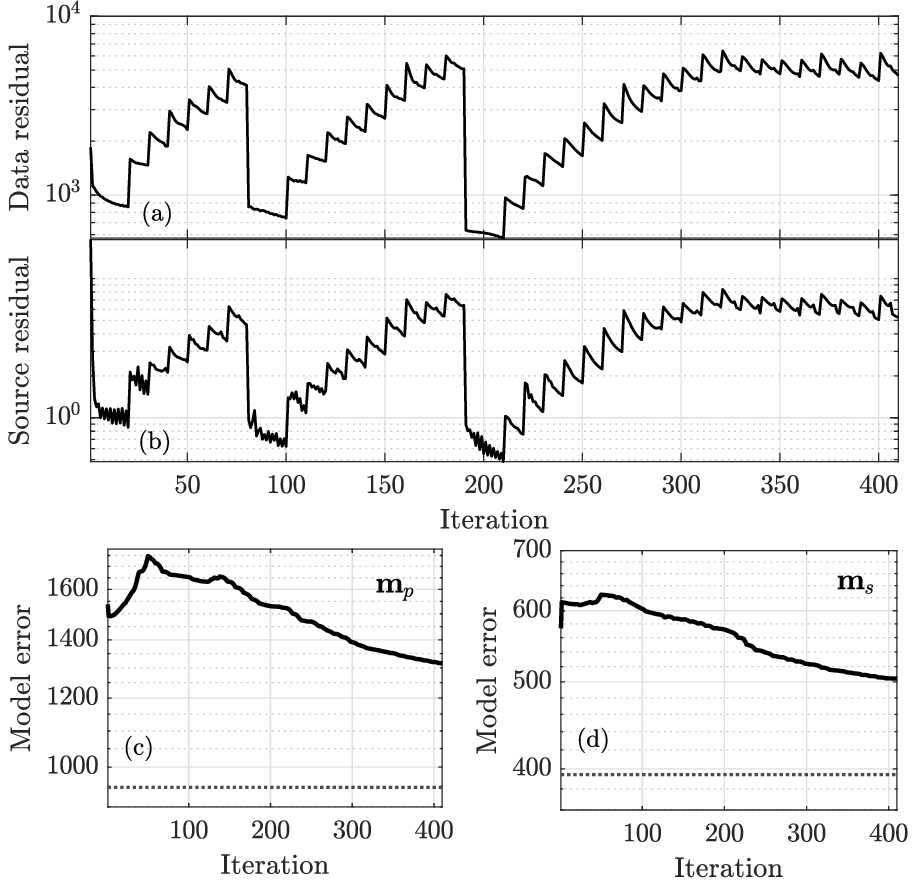}
\caption{Inversion with free surface effects. 
The variation with iteration associated to (a) data residual, (b) source residual, (c) $\bold{m}_p$ model error, and (d) $\bold{m}_s$ model error. The horizontal
dashed line represents the final model error achieved without free
surface effects.} 
\label{Fig:Over_Free_surf_error}
\end{figure}

\noindent
In this subsection, we investigate the impact of free surface effects on the accurate estimation of subsurface properties in EFWI. It is well-known that free surface effects can compromise the reliability of the inversion process \citep{brossier2009seismic}. To illustrate this, we compare computed seismograms with and without free surface effects in Figures~\ref{Fig:shot_FS}a-\ref{Fig:shot_FS}d, respectively. The presence of high-energy free surface effects introduces additional nonlinearities to the problem, making the inversion more challenging.
We initiated the inversion process with 1D initial models, as shown in Figure~\ref{Fig:over_initial}a, and the final inversion results {as well as estimated Poisson's ratio}, are displayed in Figure~\ref{Fig:inv_res_FS}a-c. Comparing these results with the ones obtained using absorbing boundary conditions on the surface (Figure~\ref{Fig:over_ALF}), it is evident that the presence of free surface effects has degraded the quality of the estimates. This degradation is further demonstrated by comparing the convergence curves associated with the inversion process, as shown in Figures \ref{Fig:Over_Free_surf_error} and \ref{Fig:over_ALF_vs_errors}.
Overall, these findings highlight the importance of accounting for free surface effects in EFWI and the challenges they pose in accurately estimating subsurface properties.

%
%
%
%

%
\newpage
%
%
\section{Conclusions}
We have introduced a new elastic full-waveform inversion
(EFWI) algorithm, leveraging the ADMM and reconstructed wavefields. This innovative approach addresses the intricate task of accurate subsurface property estimation with notable advantages. The
algorithm stands out for its inherent flexibility, accommodating
physical constraints, and efficient implementation of the Hessian
matrix through ADMM, mitigating interparameter cross-talk. A
key strength lies in its freedom from step length tuning, simplifying
the implementation process. Extensive numerical experiments
provided underscore the algorithm's effectiveness, showcasing
superior convergence properties and stability compared with the
classical FWI. Notably, it exhibited robustness against rough initial
models, noise, and free surface effects, ensuring reliable inversion
results. We have investigated the role of the Hessian matrix in suppressing cross-talk and the application of randomized source
sketching for enhanced efficiency. Although regularization was
not implemented in this study, the algorithm's inherent flexibility
allows for the straightforward incorporation of regularization techniques, empowering users with greater control over the inversion
process. 

\section{ACKNOWLEDGMENTS}
This research was partially funded by the SONATA BIS (grant
no. 2022/46/E/ST10/00266) of the National Science Center in
Poland.
\section{DATA AND MATERIALS AVAILABILITY}
Data associated with this research are available and can be
obtained by contacting the corresponding author.
\newpage

\appendix
%
%
\section*{Appendix A: Model update for ($\bold{V}_P, \bold{V}_S$)}
The optimization over ($\bold{V}_P, \bold{V}_S$) will be nonlinear thus we update the parameters by using the Gauss-Newton algorithm as
\begin{equation}\label{m_update:nonlinear}
\m^{k+1} = \m^{k}+\alpha^{k} \delta \m^{k},
\end{equation}
where $\alpha^k$ is the step-length and $\delta \m$ is the step direction, satisfying
%
\begin{equation}
\bold{H}^k\delta \m^k = -\bold{g}^k,
\end{equation}
where $\bold{g}^k$ and $\bold{H}^k$ are the gradient and the Hessian of the AL function $\mathcal{L}_{\beta}(\u^{k+1},\m,\bold{s}^{k})$ defined as:
\begin{subequations}\label{eq:Hg}
\begin{align}
&\bold{g}^{k} = [(\partial_{\bold{m}}\bold{A})\u^{k+1}]^{T}[\A\u^{k+1}-\b-\bold{s}^{k}], \label{eq:Hg_grad} \\ 
&\bold{H}^{k} = [(\partial_{\bold{m}}\bold{A}) \u^{k+1}]^{T}[(\partial_{\bold{m}}\bold{A})\u^{k+1}], \label{eq:Hg_Hess}
\end{align}
\end{subequations}
where $\A \equiv \A(\m^{k})$. Here the required term is computing $\partial_{\bold{m}}\bold{A}$, which is: 
\begin{equation}
\partial_{\bold{m}}\bold{A} = \begin{bmatrix}
\partial_{\bold{V}_\text{P}}\bold{A}_{x} & \partial_{\bold{V}_\text{S}}\bold{A}_{x} \\
\partial_{\bold{V}_\text{P}}\bold{A}_{z} & \partial_{\bold{V}_\text{S}}\bold{A}_{z} 
\end{bmatrix}
\end{equation}
where:
\begin{equation}
\begin{aligned}
& \partial_{\bold{V}_\text{P}}\bold{A}_{x} = 2\Bigl[ \diag \Bigl( \partial_{xx} \circ(\bs{\rho}\circ\bold{V}_\text{P}^k) \Bigr) +\diag \Bigl( \partial_{xz} \circ (\bs{\rho} \circ \bold{V}_\text{P}^k)\Bigr) \Bigr], \\
& \partial_{\bold{V}_\text{S}}\bold{A}_{x} = 2\Bigl[ \diag\Bigl( \partial_{zz} \circ(\bs{\rho}\circ\bold{V}_\text{S}^k) \Bigr)-\diag\Bigl( \partial_{xz} \circ(\bs{\rho}\circ\bold{V}_\text{S}^k) \Bigr) \Bigr], \\
& \partial_{\bold{V}_\text{P}}\bold{A}_{z} = 2\Bigl[\diag\Bigl( \partial_{zz} \circ(\bs{\rho}\circ\bold{V}_\text{P}^k) \Bigr) +\diag\Bigl( \partial_{xz} \circ(\bs{\rho}\circ\bold{V}_\text{P}^k) \Bigr) \Bigr], \\
& \partial_{\bold{V}_\text{S}}\bold{A}_{z} = 2\Bigl[\diag\Bigl( \partial_{xx} \circ(\bs{\rho}\circ\bold{V}_\text{S}^k) \Bigr)-\diag \Bigl( \partial_{xz} \circ(\bs{\rho}\circ\bold{V}_\text{S}^k) \Bigr) \Bigr]. \\
\end{aligned}
\end{equation} 
%
%
%
%

\section*{Appendix B: Projection onto intersection of two convex sets by ADMM}\label{Appendix:projection}
Given two closed convex sets $C_1$ and $C_2$, the model update step in equation \ref{ADMM_m} requires solving
\begin{equation}\label{m_update_cnvx}
\min_{\bold{m}}~ \|\bold{L}\bold{m}-\bold{y} \|_{2}^{2} \quad \text {s.t.} \quad \m \in C_1 \cap C_2.
\end{equation}
The constrained problem described above can be expressed as an equivalent optimization problem:
\begin{equation}\label{best_approx_prob}
\min_{\bold{m}}~ \|\bold{L}\bold{m}-\bold{y} \|_{2}^{2}+I_C{(\m)},
\end{equation}
where $C$ is the intersection of two sets $C_1$ and $C_2$, denoted by $C_1 \cap C_2$, and $I_C{(\m)}$ represent the indicator function of $C$:
\begin{equation}
I_C{(\m)}=\left\{\begin{array}{l}
0 \quad \text { if } \quad \m \in C \\
\infty\quad \text { if }\quad \m \notin C.
\end{array}\right.
\end{equation}
Including the indicator functions guarantee that any candidate solution $\m$ that is not in $C$ will have an infinite objective value, making it infeasible.
To solve \ref{best_approx_prob}, we split the two terms of the objective function by introducing auxiliary variables $\bold{p}$ and constraint $\m=\bold{p}$:
\begin{equation}\label{best_approx_probb}
\min_{\bold{m},\bold{p}}~ \|\bold{L}\bold{m}-\bold{y} \|_{2}^{2}+I_C{(\bold{p})}\quad 
 {s.t.}\quad \m-\bold{p}=0.
\end{equation}
%
%
%
%
%
%
The ADMM algorithm solves this constrained problem by the following iteration:
\begin{subequations}\label{ADMM_cnvx}
\begin{align}
\small
\m^{k+1} &= (\bold{L}^T\bold{L}+\gamma \bold{I})^{-1}(\bold{L}^T\bold{y}+\gamma(\bold{p}^{k}+\bold{q}^{k})),
 \label{ADMM_cnvx_m} \\
\bold{p}^{k+1} &= \text{P}_{C}(\m^{k+1}-\bold{q}^{k}),  \label{ADMM_cnvx_p} \\
\bold{q}^{k+1} &= \bold{q}^{k}+ \bold{p}^{k+1} -\m^{k+1}, \label{ADMM_cnvx_eta2}
\end{align}
\end{subequations} 
where $\bold{q}$ is the scaled dual variable, $\gamma>0$ is the penalty parameter, and $ \text{P}_{C}$ is the projection operators onto $C$.

\begin{thebibliography}{}
\itemsep0pt

\bibitem[Aghamiry et~al., 2019]{aghamiry2019improving}
Aghamiry, H.~S., A. Gholami, and S. Operto,  2019, Improving full-waveform
  inversion by wavefield reconstruction with the alternating direction method
  of multipliers: Geophysics, {\bf 84}, R139--R162.

\bibitem[Aghamiry et~al., 2020]{aghamiry2020multiparameter}
--------, 2020, Multiparameter wavefield reconstruction inversion for wavespeed
  and attenuation with bound constraints and total variation regularization:
  Geophysics, {\bf 85}, R381--R396.

\bibitem[Aghazade et~al., 2021]{Aghazade_2021_RSS}
Aghazade, K., H.~S. Aghamiry, A. Gholami, and S. Operto,  2021, Randomized
  source sketching for full waveform inversion: IEEE Transactions on Geoscience
  and Remote Sensing, {\bf 60}, 1--12.

\bibitem[Barnes and Charara, 2008]{Barnes_2008_FWI}
Barnes, C., and M. Charara,  2008, Full-waveform inversion results when using
  acoustic approximation instead of elastic medium: SEG Technical Program
  Expanded Abstracts, {\bf 27}, 1895--1899.

\bibitem[Barnes and Charara, 2009]{Barnes_2009_DAA}
--------, 2009, The domain of applicability of acoustic full-waveform inversion
  for marine seismic data: Geophysics, {\bf 74}, WCC91--WCC103.

\bibitem[Baumstein, 2013]{baumstein2013pocs}
Baumstein, A.,  2013, {POCS}-based geophysical constraints in multi-parameter
  full wavefield inversion: 75$^{th}$ EAGE Conference \& Exhibition
  incorporating SPE EUROPEC 2013, European Association of Geoscientists \&
  Engineers, cp--348.

\bibitem[Borisov et~al., 2022]{borisov2022graph}
Borisov, D., R.~D. Miller, S.~L. Peterie, J. Ivanov, A.~M. Hoch, and S.~D.
  Sloan,  2022, Graph-space optimal transport-based {3D} elastic {FWI} for
  near-surface seismic applications: Second International Meeting for Applied
  Geoscience \& Energy, Society of Exploration Geophysicists and American
  Association of Petroleum, 2153--2157.

\bibitem[Borisov et~al., 2014]{borisov2014acoustic}
Borisov, D., A. Stopin, and R. Plessix,  2014, Acoustic pseudo-density full
  waveform inversion in the presence of hard thin beds: 76th EAGE Conference
  and Exhibition 2014, EAGE Publications BV, 1--5.

\bibitem[Boyd et~al., 2011]{boyd2011distributed}
Boyd, S., N. Parikh, and E. Chu,  2011, Distributed optimization and
  statistical learning via the alternating direction method of multipliers: Now
  Publishers Inc.

\bibitem[Brocher, 2005]{brocher2005empirical}
Brocher, T.~M.,  2005, Empirical relations between elastic wavespeeds and
  density in the {E}arth's crust: Bulletin of the seismological Society of
  America, {\bf 95}, 2081--2092.

\bibitem[Brossier et~al., 2009]{brossier2009seismic}
Brossier, R., S. Operto, and J. Virieux,  2009, Seismic imaging of complex
  onshore structures by 2{D} elastic frequency-domain full-waveform inversion:
  Geophysics, {\bf 74}, WCC105--WCC118.

\bibitem[Bunks et~al., 1995]{bunks1995multiscale}
Bunks, C., F.~M. Saleck, S. Zaleski, and G. Chavent,  1995, Multiscale seismic
  waveform inversion: Geophysics, {\bf 60}, 1457--1473.

\bibitem[Chen et~al., 2022]{chen2022elastic}
Chen, G., W. Yang, H. Wang, H. Zhou, and X. Huang,  2022, Elastic full waveform
  inversion based on full-band seismic data reconstructed by dual
  deconvolution: IEEE Geoscience and Remote Sensing Letters, {\bf 19}, 1--5.

\bibitem[Chen and Cao, 2016]{chen2016modeling}
Chen, J.-B., and J. Cao,  2016, Modeling of frequency-domain elastic-wave
  equation with an average-derivative optimal method: Geophysics, {\bf 81},
  T339--T356.

\bibitem[Chi et~al., 2014]{chi2014full}
Chi, B., L. Dong, and Y. Liu,  2014, Full waveform inversion method using
  envelope objective function without low frequency data: Journal of Applied
  Geophysics, {\bf 109}, 36--46.

\bibitem[Choi et~al., 2008]{choi2008frequency}
Choi, Y., D.-J. Min, and C. Shin,  2008, Frequency-domain elastic full waveform
  inversion using the new pseudo-{H}essian matrix: Experience of elastic
  {M}armousi-2 synthetic data: Bulletin of the Seismological Society of
  America, {\bf 98}, 2402--2415.

\bibitem[Duan and Sava, 2016]{duan2016elastic}
Duan, Y., and P. Sava,  2016, Elastic wavefield tomography with physical model
  constraints: Geophysics, {\bf 81}, R447--R456.

\bibitem[Dykstra, 1983]{dykstra1983algorithm}
Dykstra, R.~L.,  1983, An algorithm for restricted least squares regression:
  Journal of the American Statistical Association, {\bf 78}, 837--842.

\bibitem[Forgues and Lambar\'e, 1997]{Forgues_1997_PSA}
Forgues, E., and G. Lambar\'e,  1997, Parameterization study for acoustic and
  elastic ray+{B}orn inversion: Journal of Seismic Exploration, {\bf 6},
  253--278.

\bibitem[Gabay and Mercier, 1976]{Gabay_1976_ADA}
Gabay, D., and B. Mercier,  1976, A dual algorithm for the solution of
  nonlinear variational problems via finite element approximation: Computers \&
  mathematics with applications, {\bf 2}, 17--40.

\bibitem[Gauthier et~al., 1986]{Gauthier_1986_TNI}
Gauthier, O., J. Virieux, and A. Tarantola,  1986, Two-dimensional nonlinear
  inversion of seismic waveform : numerical results: Geophysics, {\bf 51},
  1387--1403.

\bibitem[Gholami, 2023]{Gholami_2023_FWI}
Gholami, A.,  2023, Full waveform inversion and {L}agrange multipliers: arXiv
  preprint arXiv:2311.11010.

\bibitem[Gholami et~al., 2022]{gholami2022extended}
Gholami, A., H.~S. Aghamiry, and S. Operto,  2022, Extended-space full-waveform
  inversion in the time domain with the augmented {L}agrangian method:
  Geophysics, {\bf 87}, R63--R77.

\bibitem[Gholami et~al., 2023]{gholami2023multiplier}
--------, 2023, Multiplier waveform inversion ({MWI}): A reduced-space {FWI} by
  the method of multipliers: Geophysics, {\bf 88}, 1--60.

\bibitem[Gholami et~al., 2013]{Gholami_2013_WPA1}
Gholami, Y., R. Brossier, S. Operto, A. Ribodetti, and J. Virieux,  2013, Which
  parametrization is suitable for acoustic {VTI} full waveform inversion? -
  {P}art 1: sensitivity and trade-off analysis: Geophysics, {\bf 78},
  R81--R105.

\bibitem[Haber et~al., 2000]{haber2000optimization}
Haber, E., U.~M. Ascher, and D. Oldenburg,  2000, On optimization techniques
  for solving nonlinear inverse problems: Inverse problems, {\bf 16}, 1263.

\bibitem[Hu et~al., 2022]{hu2022frequency}
Hu, Y., L.-Y. Fu, Q. Li, W. Deng, and L. Han,  2022, Frequency-wavenumber
  domain elastic full waveform inversion with a multistage phase correction:
  Remote Sensing, {\bf 14}, 5916.

\bibitem[Jun et~al., 2014]{jun2014laplace}
Jun, H., Y. Kim, J. Shin, C. Shin, and D.-J. Min,  2014,
  Laplace-{F}ourier-domain elastic full-waveform inversion using time-domain
  modeling: Geophysics, {\bf 79}, R195--R208.

\bibitem[Kazei and Alkhalifah, 2019]{kazei2019scattering}
Kazei, V., and T. Alkhalifah,  2019, Scattering radiation pattern atlas: What
  anisotropic elastic properties can body waves resolve?: Journal of
  Geophysical Research: Solid Earth, {\bf 124}, 2781--2811.

\bibitem[K{\"o}hn et~al., 2012]{kohn2012influence}
K{\"o}hn, D., D. De~Nil, A. Kurzmann, A. Przebindowska, and T. Bohlen,  2012,
  On the influence of model parametrization in elastic full waveform
  tomography: Geophysical Journal International, {\bf 191}, 325--345.

\bibitem[Kwon et~al., 2017]{kwon2017waveform}
Kwon, J., H. Jun, H. Song, U.~G. Jang, and C. Shin,  2017, Waveform inversion
  in the shifted {L}aplace domain: Geophysical Journal International, {\bf 210},
  340--353.

\bibitem[Li and Demanet, 2015]{li2015phase}
Li, Y.~E., and L. Demanet,  2015, Phase and amplitude tracking for seismic
  event separation: Geophysics, {\bf 80}, WD59--WD72.

\bibitem[Lin and Huang, 2014]{lin2014acoustic}
Lin, Y., and L. Huang,  2014, Acoustic-and elastic-waveform inversion using a
  modified total-variation regularization scheme: Geophysical Journal
  International, {\bf 200}, 489--502.

\bibitem[Marty et~al., 2022]{marty2022elastic}
Marty, P., C. Boehm, and A. Fichtner,  2022, Elastic full-waveform inversion
  for transcranial ultrasound computed tomography using optimal transport:
  IEEE International Ultrasonics Symposium (IUS), IEEE, 1--4.

\bibitem[M{\'e}tivier et~al., 2014]{Metivier_2014_MFW}
M{\'e}tivier, L., R. Brossier, S. Operto, and J. Virieux,  2014,
  Multi-parameter {FWI}-an illustration of the {H}essian operator role for
  mitigating trade-off between parameter classes: 6$^{th}$ EAGE Saint
  Petersburg International Conference and Exhibition, 1--5.

\bibitem[Mora, 1988]{Mora_1988_EWI}
Mora, P.~R.,  1988, Elastic wavefield inversion of reflection and transmission
  data: Geophysics, {\bf 53}, 750--759.

\bibitem[Mulder and Plessix, 2008]{Mulder_2008_ESI}
Mulder, W., and R.~E. Plessix,  2008, Exploring some issues in acoustic full
  waveform inversion: Geophysical Prospecting, {\bf 56}, 827--841.

\bibitem[Nocedal and Wright, 2006]{nocedal2006numerical}
Nocedal, J., and S. Wright,  2006, Numerical optimization: Springer Science \&
  Business Media.

\bibitem[Operto et~al., 2023]{Operto_2023_FWI}
Operto, S., A. Gholami, H.~S. Aghamiry, G. Guo, F. Mamfoumbi, S. Beller, K.
  Aghazade, F. Mamfoumbi, L. Combe, and A. Ribodetti,  2023, Extending the
  search space of full-waveform inversion beyond the single-scattering born
  approximation: A tutorial review: Geophysics, {\bf 88}, 1--32.

\bibitem[Operto et~al., 2013]{Operto_2013_TLE}
Operto, S., Y. Gholami, V. Prieux, A. Ribodetti, R. Brossier, L. Metivier, and
  J. Virieux,  2013, A guided tour of multiparameter full-waveform inversion
  with multicomponent data: From theory to practice: The leading edge, {\bf
  32}, 1040--1054.

\bibitem[Operto et~al., 2006]{Operto_2006_CIM}
Operto, S., J. Virieux, J.~X. Dessa, and G. Pascal,  2006, Crustal imaging from
  multifold ocean bottom seismometers data by frequency-domain full-waveform
  tomography: application to the eastern {N}ankai trough: Journal of
  Geophysical Research, {\bf 111}.

\bibitem[Pan et~al., 2018]{pan2018elastic}
Pan, W., K.~A. Innanen, and Y. Geng,  2018, Elastic full-waveform inversion and
  parametrization analysis applied to walk-away vertical seismic profile data
  for unconventional (heavy oil) reservoir characterization: Geophysical
  Journal International, {\bf 213}, 1934--1968.

\bibitem[Pan et~al., 2019]{pan2019interparameter}
Pan, W., K.~A. Innanen, Y. Geng, and J. Li,  2019, Interparameter trade-off
  quantification for isotropic-elastic full-waveform inversion with various
  model parameterizations: Geophysics, {\bf 84}, R185--R206.

\bibitem[Powell, 1969]{Powell_1969_NLC}
Powell, M.~J.,  1969, A method for nonlinear constraints in minimization
  problems: Optimization,  283--298.

\bibitem[Prieux et~al., 2011]{Prieux_2011_FAI}
Prieux, V., R. Brossier, Y. Gholami, S. Operto, J. Virieux, O. Barkved, and J.
  Kommedal,  2011, On the footprint of anisotropy on isotropic full waveform
  inversion: the {V}alhall case study: Geophysical Journal International, {\bf
  187}, 1495--1515.

\bibitem[Prieux et~al., 2013a]{Prieux_2013_MFWa}
Prieux, V., R. Brossier, S. Operto, and J. Virieux,  2013a, Multiparameter full
  waveform inversion of multicomponent {OBC} data from {V}alhall. {P}art 1:
  imaging compressional wavespeed, density and attenuation: Geophysical Journal
  International, {\bf 194}, 1640--1664.

\bibitem[Prieux et~al., 2013b]{Prieux_2013_MFWb}
--------, 2013b, Multiparameter full waveform inversion of multicomponent {OBC}
  data from valhall. {P}art 2: imaging compressional and shear-wave velocities:
  Geophysical Journal International, {\bf 194}, 1665--1681.

\bibitem[Ravaut et~al., 2004]{Ravaut_2004_MSI}
Ravaut, C., S. Operto, L. Improta, J. Virieux, A. Herrero, and P.
  dell'Aversana,  2004, Multi-scale imaging of complex structures from
  multi-fold wide-aperture seismic data by frequency-domain full-wavefield
  inversions: application to a thrust belt: Geophysical Journal International,
  {\bf 159}, 1032--1056.

\bibitem[Sears et~al., 2008]{Sears_2008_EFW}
Sears, T., S. Singh, and P. Barton,  2008, {Elastic full waveform inversion of
  multi-component OBC seismic data}: Geophysical Prospecting, {\bf 56},
  843--862.

\bibitem[Sun and Demanet, 2021]{sun2021deep}
Sun, H., and L. Demanet,  2021, Deep learning for low-frequency extrapolation
  of multicomponent data in elastic {FWI}: IEEE Transactions on Geoscience and
  Remote Sensing, {\bf 60}, 1--11.

\bibitem[Sun et~al., 2017]{sun2017density}
Sun, M., J. Yang, L. Dong, Y. Liu, and C. Huang,  2017, Density reconstruction
  in multiparameter elastic full-waveform inversion: Journal of Geophysics and
  Engineering, {\bf 14}, 1445--1462.

\bibitem[Tarantola, 1984]{Tarantola_1984_ISR}
Tarantola, A.,  1984, Inversion of seismic reflection data in the acoustic
  approximation: Geophysics, {\bf 49}, 1259--1266.

\bibitem[Tarantola, 1986]{Tarantola_1986_SNL}
--------, 1986, A strategy for nonlinear inversion of seismic reflection data:
  Geophysics, {\bf 51}, 1893--1903.

\bibitem[{van Leeuwen} and Herrmann, 2013]{VanLeeuwen_2013_MLM}
{van Leeuwen}, T., and F.~J. Herrmann,  2013, Mitigating local minima in
  full-waveform inversion by expanding the search space: Geophysical Journal
  International, {\bf 195}, 661--667.

\bibitem[Vigh et~al., 2014]{vigh2014elastic}
Vigh, D., K. Jiao, D. Watts, and D. Sun,  2014, Elastic full-waveform inversion
  application using multicomponent measurements of seismic data collection:
  Geophysics, {\bf 79}, R63--R77.

\bibitem[Virieux and Operto, 2009]{Virieux_2009_OFW}
Virieux, J., and S. Operto,  2009, An overview of full waveform inversion in
  exploration geophysics: Geophysics, {\bf 74}, WCC1--WCC26.
\bibitem[Wang et~al., 2019]{wang2019elastic}
Wang, C., J. Cheng, W.~W. Weibull, and B. Arntsen,  2019, Elastic wave-equation
  migration velocity analysis preconditioned through mode decoupling:
  Geophysics, {\bf 84}, R341--R353.
  \bibitem[Wang and Cheng, 2017]{wang2017elastic}
Wang, T., and J. Cheng,  2017, Elastic full waveform inversion based on mode
  decomposition: The approach and mechanism: Geophysical Journal International,
  {\bf 209}, 606--622.
  \bibitem[Wang et~al., 2018]{wang2018elastic}
Wang, T., J. Cheng, Q. Guo, and C. Wang,  2018, Elastic wave-equation-based
  reflection kernel analysis and traveltime inversion using wave mode
  decomposition: Geophysical Journal International, {\bf 215}, 450--470.
\bibitem[Xiong et~al., 2023]{xiong2023improved}
Xiong, K., D. Lumley, and W. Zhou,  2023, Improved seismic envelope full
  waveform inversion: Geophysics, {\bf 88}, 1--197.

\bibitem[Xu and McMechan, 2014]{xu20142d}
Xu, K., and G.~A. McMechan,  2014, 2{D} frequency-domain elastic full-waveform
  inversion using time-domain modeling and a multistep-length gradient
  approach: Geophysics, {\bf 79}, R41--R53.

\bibitem[Zhang and Li, 2022]{zhang2022elastic}
Zhang, W., and Y. Li,  2022, Elastic wave full-waveform inversion in the time
  domain by the trust region method: Journal of Applied Geophysics, {\bf 197},
  104540.

\bibitem[Zhang and Alkhalifah, 2019]{zhang2019local}
Zhang, Z.-d., and T. Alkhalifah,  2019, Local-crosscorrelation elastic
  full-waveform inversion: Geophysics, {\bf 84}, R897--R908.

\bibitem[Zhang et~al., 2018]{zhang2018multiparameter}
Zhang, Z.-d., T. Alkhalifah, E.~Z. Naeini, and B. Sun,  2018, Multiparameter
  elastic full waveform inversion with facies-based constraints: Geophysical
  Journal International, {\bf 213}, 2112--2127.

\end{thebibliography}

\end{document}